\documentclass[11pt]{amsart}
\usepackage{epsfig,amsthm}
\usepackage{pictex}

\DeclareFontFamily{OMX}{cmex}{}{}
\DeclareFontShape{OMX}{cmex}{m}{n}{
  <6> <7> <8> <10> <10.95> <12> <14.4> <17.5> cmex10
  }{}
\makeatletter
\outer\def\proclaim #1. #2\par{%
  \medbreak
  \noindent{\bfseries#1.\enspace}{\slshape#2\par}%
  \ifdim\lastskip<\medskipamount
    \removelastskip\penalty55\medskip
  \fi}

{\theoremstyle{definition}
\thm@bodyfont{\rmfamily}
\newtheorem{Rem}{Remark}
\newtheorem{Def}{Definition}
\newtheorem{Not}{Notation}
\newtheorem{Ex}{Example}}

\newtheorem{Lem}{Lemma}
\newtheorem{Prop}{Proposition}
\newtheorem{Thm}{Theorem}
\newtheorem{Tad}{Theorem-and-Definition}
\newtheorem{Cor}{Corollary}
\makeatother

\newcommand{\capt}{\bf T \rm}  
\newcommand{\tzero}{\bf T_0\rm}
\newcommand{\capa}{\bf A \rm}  

\newcommand{\rootspace}{\bf\Phi\rm} 
\newcommand{\capz}{\bf Z \rm}  
\newcommand{\capq}{\bf Q \rm}  
\newcommand{\capc}{\bf C \rm}  
\newcommand{\capr}{\bf R \rm}  
\newcommand{\ssubp}{\bf S_p \rm}  
\newcommand{\ssubfour}{\bf S_4 \rm}  

\begin{document}
\title{An Obstruction to Embedding 4-tangles in Links} 

\author{David A. Krebes} 

\thanks{This paper was prepared while the author was a graduate
student at the University of Illinois at Chicago, a Pacific Institute
for the Mathematical Sciences post-doctoral fellow at the University
of British Columbia, and a visiting assistant professor at the
University of South Alabama.}

\begin{abstract}
We consider the ways in which a $4$-tangle $T$ inside a unit cube can
be extended outside the cube into a knot or link $L$.  We present two
links $n(T)$ and $d(T)$ such that the greatest common divisor of the
determinants of these two links always divides the determinant of the
link $L$.

In order to prove this result we give a two-integer invariant of
$4$-tangles.  Calculations are facilitated by viewing the determinant
as the Kauffman bracket at a fourth root of -1, which sets the loop
factor to zero.  For rational tangles, our invariant coincides with
the value of the associated continued fraction.

\end{abstract}

\maketitle

\parindent = 0.5 in

\section{Introduction}
\label{sec:introduction}

In this article, we will consider the ways in which a knot or link $L$
can intersect a ball $B$ in $\bf R^3$.  We examine the case in which
$L$ meets the boundary of $B$ transversely in four points. The pair
$(B,L\cap B)$ is called a {\em $4$-tangle} (or simply {\em
tangle}) and is said to {\em sit inside} or be {\em embedded in} the
link $L$. Section~\ref{sec:domain} describes tangles.  (Note that the
definition of tangle which we will use here is slightly more general
than that found in \cite{conway}, which deals only with what we will
call loop-free tangles.)

\begin{figure}
\begin{center}
\font\thinlinefont=cmr5
\begingroup\makeatletter\ifx\SetFigFont\undefined%
\gdef\SetFigFont#1#2#3#4#5{%
  \reset@font\fontsize{#1}{#2pt}%
  \fontfamily{#3}\fontseries{#4}\fontshape{#5}%
  \selectfont}%
\fi\endgroup%
\mbox{\beginpicture
\setcoordinatesystem units <1.00000cm,1.00000cm>
\unitlength=1.00000cm
\linethickness=1pt
\setplotsymbol ({\makebox(0,0)[l]{\tencirc\symbol{'160}}})
\setshadesymbol ({\thinlinefont .})
\setlinear
%
%
\linethickness= 0.500pt
\setplotsymbol ({\thinlinefont .})
\putrectangle corners at  2.381 26.194 and  3.969 25.082
%
%
\linethickness= 0.500pt
\setplotsymbol ({\thinlinefont .})
\putrectangle corners at  8.572 26.194 and 10.160 25.082
\linethickness= 0.500pt
\setplotsymbol ({\thinlinefont .})
%
%
%
\plot    2.381 25.876  2.302 25.876
         2.223 25.896
         2.143 25.956
         2.103 26.075
         2.113 26.164
         2.143 26.273
         2.237 26.377
         2.325 26.418
         2.441 26.452
         2.509 26.466
         2.583 26.478
         2.665 26.488
         2.753 26.496
         2.849 26.503
         2.951 26.508
         3.059 26.510
         3.175 26.511
         3.291 26.511
         3.399 26.509
         3.501 26.506
         3.597 26.501
         3.685 26.496
         3.767 26.489
         3.841 26.481
         3.909 26.472
         4.025 26.449
         4.113 26.422
         4.207 26.352
         4.242 26.278
         4.266 26.214
         4.286 26.114
         4.247 26.035
         4.128 25.956
         /
\plot  4.128 25.956  3.969 25.876 /
\linethickness= 0.500pt
\setplotsymbol ({\thinlinefont .})
%
%
%
\plot    2.381 25.400  2.302 25.400
         2.223 25.360
         2.143 25.241
         2.113 25.167
         2.103 25.102
         2.143 25.003
         2.242 24.924
         2.307 24.884
         2.381 24.844
         2.490 24.810
         2.567 24.796
         2.659 24.785
         2.766 24.776
         2.887 24.770
         2.954 24.768
         3.024 24.766
         3.097 24.765
         3.175 24.765
         3.253 24.765
         3.328 24.766
         3.399 24.768
         3.468 24.770
         3.533 24.773
         3.595 24.776
         3.711 24.785
         3.814 24.796
         3.904 24.810
         3.982 24.826
         4.048 24.844
         4.152 24.889
         4.227 24.944
         4.271 25.008
         4.286 25.082
         4.276 25.157
         4.247 25.221
         4.128 25.321
         /
\plot  4.128 25.321  3.969 25.400 /
\linethickness= 0.500pt
\setplotsymbol ({\thinlinefont .})
%
%
%
\plot    8.572 25.876  8.493 25.876
         8.414 25.856
         8.334 25.797
         8.295 25.698
         8.305 25.633
         8.334 25.559
         8.414 25.440
         8.493 25.400
         /
\plot  8.493 25.400  8.572 25.400 /
\linethickness= 0.500pt
\setplotsymbol ({\thinlinefont .})
%
%
%
\plot   10.160 25.876 10.319 25.797
        10.438 25.718
        10.478 25.638
        10.438 25.559
        10.319 25.479
         /
\plot 10.319 25.479 10.160 25.400 /
%
%
\put{\SetFigFont{17}{20.4}{\rmdefault}{\mddefault}{\updefault}$n(T)=$} [lB] at  0.1 25.559
%
%
\put{\SetFigFont{17}{20.4}{\rmdefault}{\mddefault}{\updefault}$d(T)=$} [lB] at  6.23 25.559
%
%
\put{\SetFigFont{17}{20.4}{\rmdefault}{\mddefault}{\updefault};} [lB] at  4.921 25.559
%
%
\put{\SetFigFont{20}{24.0}{\rmdefault}{\mddefault}{\updefault}T} [lB] at  9.049 25.400
%
%
\put{\SetFigFont{20}{24.0}{\rmdefault}{\mddefault}{\updefault}T} [lB] at  2.857 25.400
\linethickness=0pt
\putrectangle corners at  0.159 26.528 and 10.494 24.748
\endpicture}
\end{center}
\caption{Two important links in which the tangle $T$ embeds}
\label{fig:numden_links}
\end{figure}

The results proven in this paper involve the {\em determinant} of a
knot or link $L$, which is a non-negative integer denoted $|\langle
L\rangle |$.  (This notation has a meaning which will be explained
later in this introduction; see also Section~\ref{sec:determinant}.)

We are interested in the ways a tangle $T$ can be
embedded in a link.  Two such embeddings are given in
Figure~\ref{fig:numden_links}.  

The main result of our paper is the following:

\proclaim{Theorem \ref{Thm:9main} (Version 2)}.  If the tangle $T$ can be embedded in the link $L$, then the greatest common divisor of the determinants of $n(T)$ and of $d(T)$ divides the determinant of $L$.

Since the determinant of the unknot is $1$, we have the following corollary.

\proclaim{Corollary \ref{Cor:unknot}}.  The tangle $T$ cannot sit
inside the unknot unless the determinants of $n(T)$ and $d(T)$ are
relatively prime.

\begin{Rem}
Corollary \ref{Cor:unknot} can also be regarded as a knotting
criterion: If $K\cong S^1\subset S^3$ has a tangle $T$ embedded in it
such that the determinants of $n(T)$ and $d(T)$ are not relatively
prime, then $K$ is knotted.
\end{Rem}

\begin{figure}
\begin{center}
\leavevmode
\epsfbox{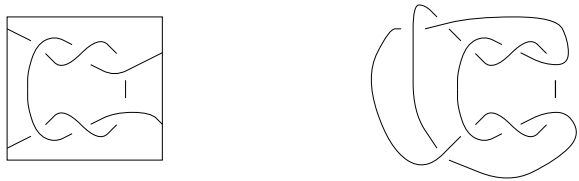}
\end{center} 
\caption{The square tangle and a knot in which it sits}
\label{fig:intro}
\end{figure}

The following claim illustrates Corollary \ref{Cor:unknot}: If a closed curve
intersects the unit cube (a topological ball) in the tangle $T$ shown
on the left in Figure~\ref{fig:intro} ($T$ is known to laypeople as a square
knot), then it is genuinely knotted.

For example the curve depicted on the right in Figure~\ref{fig:intro}
is knotted: A rubber band, for instance, cannot be manipulated into
this shape without breaking it and gluing the ends back together.  The
determinant of $n(T)$, a connected sum of two trefoils, is $9$ and the
determinant of $d(T)$, the unlink of two components, is $0$.  Since
$9$ and $0$ are not relatively prime, Corollary~\ref{Cor:unknot} may
be applied, proving the claims.

The particulars of this example are contained in
Section~\ref{sec:calculations}.  An alternative proof for the reader
familiar with Fox $n$-colourings is found in
Section~\ref{sec:colourings}.

To prove Theorem \ref{Thm:9main}, Version 2, we introduce an invariant
of tangles whose values are formal fractions $p/q$, not necessarily
reduced, as defined in Section
\ref{sec:range}.  A precise definition of the invariant appears in
Section~\ref{sec:definition}.  However, it is rather easy to define $p$ and $q$ up to sign:

\proclaim{Proposition \ref{Prop:absvalue}}.  If $f(T)=p/q$ then $|p|$ and $|q|$ are the determinants $|\langle n(T)\rangle
|$ and $|\langle d(T)\rangle |$, respectively.

Thus Theorem \ref{Thm:9main}, Version 2 is equivalent to:

\proclaim{Theorem \ref{Thm:9main}, Version 1}.  If $T$ is a tangle with invariant
$p/q$ and $T$ can be embedded in the link $L$, then $gcd(p,q)$ divides
${|\langle L\rangle |}.$

To define this invariant we consider the Kauffman bracket, described
in \cite{state}, evaluated at an eighth root of unity, namely
$A=e^{\pi i/4}$ where $A$ is the indeterminate in the Kauffman
bracket polynomial.  Some consequences of this choice are described in
Section~\ref{sec:value}.  With this choice of $A$ the absolute value
$|\langle L \rangle |$ of the Kauffman bracket is precisely the
determinant of the knot or link $L$.

Our invariant has an agreeable additivity property:

\proclaim{Proposition \ref{Prop:8additive}}. (paraphrase). 
If $p/q$ is the value of the invariant on the tangle $T$, and
$r/s$ is the value of the invariant on the tangle $T'$,  then
the value of the invariant on the tangle $T+T'$, which is described in
Section~\ref{sec:domain}, is $$\frac{ps+qr}{qs}.$$

This last expression is seen to be the formal sum of $p/q$ and $r/s$ when the product of the denominators is
used as a common denominator in elementary school addition of
fractions.

Proposition \ref{Prop:8additive} is used in conjunction with the following:

\proclaim{Lemma \ref{Lem:9ambient}}.  If the tangle $T$ sits inside the
link $L$, then there is a second tangle $T'$ such that $L$ is ambient
isotopic to the link $n(T+T')$.

Thus if $p/q$ is the value of the invariant on the tangle $T$
then there are integers $r$ and $s$ such that the determinant of the
link $L$ is $|ps+qr|$.  Thus any common divisor of $p$ and $q$ divides
the determinant of $L$, and Theorem~\ref{Thm:9main} is established.

In Section \ref{sec:realiz} we illustrate the power of the invariant
by showing that any fraction $p/q$ (not necessarily reduced) is
associated to at least one tangle.  If $p$ or $q$ is odd, this tangle
may be chosen to consist of two unknotted arcs in $B$.  (It is
loop-free by Proposition~\ref{Prop:pqparity}.)

A by-product of this investigation is an efficient way of calculating
the determinants of a group of links systematically using the Kauffman
bracket.  See Section~\ref{sec:determinant}.

\section{$4$-tangles}
\label{sec:domain}

All our work will be done in the PL (piecewise-linear) category.

Consider a unit cube $J=I^3$ (here $I$ is the unit interval $[0,1]$ of
the real line) and a one-dimensional properly embedded submanifold $T$
such that $\partial T = T \cap \partial J$ consists of four points.
Clearly $T$ is homeomorphic to the disjoint union $I \coprod I \coprod
S^1 \coprod \ldots \coprod S^1$ of two unit intervals and some
number (possibly zero) of circles.  There is an isotopy of $(J,T)$ which takes
$\partial T$ to the four points
$$(0,\frac{1}{3},\frac{1}{2}),\quad (0,\frac{2}{3},\frac{1}{2}), \quad
(1,\frac{1}{3}, \frac{1}{2}),\quad (1,\frac{2}{3},\frac{1}{2})$$
of $\partial J$.  Let us therefore assume that $\partial T$ consists
of these four points.  

In this paper $J$ will be drawn as projected along the z-axis to the
square $I \times I \times \{1/2\}.$ See diagrams.

Since $T$ is homeomorphically a disjoint union it is a disconnected
space and we may refer to its components; the two components
homeomorphic to the unit interval $I$ will be called the {\em arc
components} and the components homeomorphic to $S^1$ will be called
{\em loops} or sometimes {\em cycles}.  When there are no loops the
two arc components will be referred to simply as {\em strands}.

\begin{figure}
\begin{center}
\leavevmode
\epsfbox{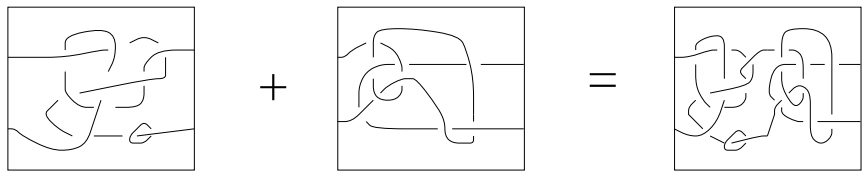}
\end{center} 
\caption{The addition rule for tangles}
\label{fig:addition}
\end{figure}

The set of such pairs $(J,T)$ has a monoid structure, where addition
is obtained by horizontal concatenation and compression, eg. see
Figure~\ref{fig:addition}.

Formally, $(J,T) + (J,T') = (J,\Lambda(T)\cup\Psi(T'))$ where
$$\Lambda: J \rightarrow J: (x,y,z) \rightarrow (\frac{1}{2} x, y,z)$$
$$\Psi: J \rightarrow J: (x,y,z) \rightarrow (\frac{1}{2} x + \frac{1}{2}, y,z)$$

Finally, a {\em $4$-tangle (tangle)} will be an ambient isotopy (deformations without
double points and with fixed endpoints) class of such pairs $(J,T)$.
The above-mentioned monoid structure on the space of pairs induces a
monoid structure on the space $\capt$ of $4$-tangles which will be
denoted ``+'', and the tangle $T+T'$ will be called the {\em sum} of
$T$ and $T'$.  The reader should recognize, however, that
this monoid structure is not commutative.

\section{Diagrams, Shadows, and States}
\label{sec:states}

Define a tangle shadow, like a link shadow (see \cite{physics}) as the
underlying combinatorial object obtained from a tangle diagram by
disregarding which is the over- and which is the undercrossing arc at
each crossing.  Link shadows can be thought of as planar graphs, where
each vertex is incident to four edges.

By a component of a tangle shadow we shall mean the part corresponding
to a component of any tangle which projects to the given shadow.  In
other words, a component is a minimal non-empty collection of edges
where the edges opposite each other at a vertex are considered to
belong to the same component.

As in the literature (eg. \cite{state}) we will use the term {\em
state} to refer both to a smoothing at a particular crossing and to a
choice of such smoothings over all crossings of a tangle or link.  If we
choose a state for each crossing of a tangle or link diagram then we
arrive at a new tangle or link diagram-- one without crossings-- which
in turn represents a (rather trivial) tangle or link.  Thus we may
also talk about components of tangle or link states.  With this in
mind we make the following definition.

\begin{Def}
Define a tangle, tangle diagram, tangle shadow or tangle state to be
{\em loop-free} or {\em acyclic} if it has only the two arc components and no
loops (like the tangles in \cite{conway}).  Two acyclic tangles,
tangle diagrams, tangle shadows or tangle states will be said to be
{\em homotopic} to one another if the strands pair off the endpoints
in the same way; equivalently, if when one starts at the upper left
hand endpoint and follows the strand originating from there one
arrives at the same endpoint in either case.
\label{Def:acyclic}
\end{Def}

\section{The Kauffman Bracket at an Eighth Root of Unity}
\label{sec:value}

See \cite{state} for an exposition of the Kauffman bracket, regular
isotopy and the Reidemeister moves.

\begin{Not} 
  Throughout this paper the expression $\langle L \rangle$
  will always mean the Kauffman bracket for the regular isotopy class
  $L$ of link diagrams {\em evaluated at}
$$A=e^{\pi i/4}.$$
\label{Not:value}
\end{Not}

\begin{Not} Let $\rootspace = \{u|u^8=1\} \subset \capc$, the
  8th roots of unity, or powers of $A$.
\end{Not}

\begin{Def}  
A {\em monocyclic} state for a link or a link shadow is a state which
  consists of a single loop.  Monocyclic states are the only states of
  interest to us since they are the only states on which the Kauffman
  bracket evaluated at $A=e^{{\pi i}/4}$ is non-zero.  This is because
  second, third or subsequent loops each introduce the factor
  $-A^2-A^{-2}$ which evaluates to $-i-(i^{-1})=-i--i=0.$
\end{Def}

\begin{Lem}  
  Let $L$ and $L'$ be two link diagrams which differ by a single
  Reidemeister move.  If this move is of type II or III then $\langle
  L\rangle=\langle L'\rangle$.  If the move is of type I then $\langle
  L\rangle =A^{\pm 1}\langle L'\rangle .$
\label{Lem:value}
\end{Lem}

\bf Proof. \rm  Invariance under type II or III moves is the meaning of the
statement that the Kauffman bracket is an invariant of regular
isotopy (see \cite{state}).

Behaviour under one version of the Type I Reidemeister move is given
in equation (\ref{eq:Rm_1}); the calculations for the other versions
(the kink may protrude to the left for instance or the roles of the
over- and undercrossing arcs may be reversed from that shown on the
left hand side of equation (\ref{eq:Rm_1})) are similar.

\begin{equation}
  \begin{split}
  \left\langle\;\;\,
  \raisebox{-0.58cm}{\mbox{%
  \beginpicture
  \setcoordinatesystem units <0.4cm,0.4cm>
  \setplotarea x from 2.1 to 4.1, y from 0 to 3.6
  \setquadratic
  \plot  2.381 1.654  2.302 1.574
         2.242 1.475
         2.227 1.411
         2.223 1.336
         2.223 1.232
         2.223 1.161
         2.223 1.078
         2.223 0.983
         2.223 0.875
         2.223 0.755
         2.223 0.690
         2.223 0.622
         2.223 0.600
         2.223 0.20
         /
  \plot  2.223 3.25  2.223 3.162
         2.223 3.064
         2.223 2.968
         2.223 2.875
         2.223 2.785
         2.223 2.697
         2.223 2.611
         2.223 2.528
         2.223 2.448
         2.225 2.369
         2.232 2.294
         2.245 2.221
         2.262 2.150
         2.285 2.082
         2.312 2.016
         2.381 1.892
         2.466 1.778
         2.560 1.674
         2.664 1.579
         2.778 1.495
         2.897 1.426
         3.016 1.376
         3.135 1.346
         3.254 1.336
         3.373 1.346
         3.493 1.376
         3.612 1.426
         3.731 1.495
         3.835 1.579
         3.909 1.674
         3.954 1.778
         3.969 1.892
         3.954 2.001
         3.909 2.090
         3.835 2.160
         3.731 2.209
         3.612 2.244
         3.493 2.269
         3.373 2.284
         3.254 2.289
         3.140 2.279
         3.036 2.249
         2.942 2.199
         2.857 2.130
         2.767 2.07
         2.6 1.91
         /
        \endpicture}}
        \;\right\rangle  
  & = A\left\langle\;\;\,
  \raisebox{-0.58cm}{\mbox{%
  \beginpicture
  \setcoordinatesystem units <0.4cm,0.4cm>
  \setplotarea x from 2.1 to 4.1, y from 0.15 to 3.75
  \setquadratic
  \plot  2.223  3.25
         2.223  3.686  2.223  3.210
         2.223  3.095
         2.223  2.987
         2.223  2.886
         2.223  2.793
         2.223  2.708
         2.223  2.630
         2.223  2.559
         2.223  2.496
         2.232  2.386
         2.262  2.297
         2.381  2.178
         2.461  2.148
         2.540  2.138
         2.619  2.148
         2.699  2.178
         2.783  2.218
         2.877  2.257
         2.982  2.297
         3.096  2.337
         3.215  2.367
         3.334  2.376
         3.453  2.367
         3.572  2.337
         3.681  2.292
         3.770  2.238
         3.840  2.173
         3.889  2.099
         3.919  2.019
         3.929  1.940
         3.919  1.861
         3.889  1.781
         3.840  1.707
         3.770  1.642
         3.681  1.588
         3.572  1.543
         3.458  1.513
         3.354  1.503
         3.259  1.513
         3.175  1.543
         3.091  1.583
         2.996  1.622
         2.892  1.662
         2.778  1.702
         2.664  1.732
         2.560  1.741
         2.466  1.732
         2.381  1.702
         2.312  1.617
         2.285  1.542
         2.262  1.444
         2.245  1.324
         2.238  1.255
         2.232  1.181
         2.228  1.101
         2.225  1.016
         2.223  0.925
         2.223  0.829
         2.223  0.2
         /
        \endpicture}}
  \;\right\rangle + A^{-1}
  \left\langle\;\;\,
  \raisebox{-0.58cm}{\mbox{%
  \beginpicture
  \setcoordinatesystem units <0.4cm,0.4cm>
  \setplotarea x from 2.1 to 4.1, y from 0 to 3.6
  \setquadratic
  \plot  2.223  3.25
         2.223  3.399  2.223  2.764
         2.223  2.686
         2.224  2.612
         2.225  2.542
         2.227  2.476
         2.234  2.354
         2.242  2.248
         2.254  2.156
         2.267  2.079
         2.302  1.970
         2.337  1.891
         2.361  1.811
         2.376  1.732
         2.381  1.652
         2.376  1.573
         2.361  1.494
         2.337  1.414
         2.302  1.335
         2.267  1.221
         2.254  1.138
         2.242  1.037
         2.234  0.920
         2.230  0.854
         2.227  0.784
         2.225  0.710
         2.224  0.632
         2.223  0.549
         2.223  0.462
         2.223  0.2
         /
%
%
\linethickness= 0.500pt
\setplotsymbol ({\tiny .})
\ellipticalarc axes ratio  0.635:0.476  360 degrees 
        from  3.969  1.652 center at  3.334  1.652
\linethickness= 0.500pt
\setplotsymbol ({\tiny .})
\endpicture}}
\;\right\rangle\\
  & = A \left\langle\;\;\,
  \raisebox{-0.58cm}{\mbox{%
  \beginpicture
  \setcoordinatesystem units <0.4cm,0.4cm>
  \setplotarea x from 2.1 to 4.1, y from 0.15 to 3.75
  \setquadratic
  \plot  2.223  3.25
         2.223  3.686  2.223  3.210
         2.223  3.095
         2.223  2.987
         2.223  2.886
         2.223  2.793
         2.223  2.708
         2.223  2.630
         2.223  2.559
         2.223  2.496
         2.232  2.386
         2.262  2.297
         2.381  2.178
         2.461  2.148
         2.540  2.138
         2.619  2.148
         2.699  2.178
         2.783  2.218
         2.877  2.257
         2.982  2.297
         3.096  2.337
         3.215  2.367
         3.334  2.376
         3.453  2.367
         3.572  2.337
         3.681  2.292
         3.770  2.238
         3.840  2.173
         3.889  2.099
         3.919  2.019
         3.929  1.940
         3.919  1.861
         3.889  1.781
         3.840  1.707
         3.770  1.642
         3.681  1.588
         3.572  1.543
         3.458  1.513
         3.354  1.503
         3.259  1.513
         3.175  1.543
         3.091  1.583
         2.996  1.622
         2.892  1.662
         2.778  1.702
         2.664  1.732
         2.560  1.741
         2.466  1.732
         2.381  1.702
         2.312  1.617
         2.285  1.542
         2.262  1.444
         2.245  1.324
         2.238  1.255
         2.232  1.181
         2.228  1.101
         2.225  1.016
         2.223  0.925
         2.223  0.829
         2.223  0.2
         /
  \endpicture}}
  \;\right\rangle\\
  & = A   \left\langle\;\;\,
  \raisebox{-0.58cm}{\mbox{%
  \beginpicture
  \setcoordinatesystem units <0.4cm,0.4cm>
  \setplotarea x from 2.1 to 4.1, y from 0 to 3.6
  \setquadratic
  \plot  2.223  3.25
         2.223  3.399  2.223  2.764
         2.223  2.686
         2.223  0.549
         2.223  0.462
         2.223  0.2
         /
   \endpicture}}
   \;\right\rangle
   \end{split}
   \label{eq:Rm_1}
   \end{equation}

Equation (\ref{eq:Rm_1}) completes the proof of Lemma~\ref{Lem:value}.

\begin{Lem}
Any two monocyclic states for a link diagram
$L$ (or for a link shadow) differ at an even number of crossings of the original diagram.
\label{Lem:even}
\end{Lem}

\bf Proof. \rm  Consider the checkerboard shading (see \cite{physics}) for
the link diagram $L$, thought of as drawn on $S^2$.  This divides $S^2$ into
shaded and unshaded regions, which we can think of as the vertices of
a graph.  A state for $L$ connects these regions in various ways,
ie. forms edges between the vertices.  The resulting planar graph is
the union of two subgraphs: the ``shaded'' one and the ``unshaded''
one.  Note that a loop of one of these graphs bounds, on either side,
at least one component of the other graph.

Now any component of either graph is bounded by a cycle of the state.
So monocyclic states induce graphs with only two components: the
shaded one and the unshaded one.  It follows that these graphs are
both trees.  Hence if there were $k$ shaded regions of the shadow,
there must be exactly $k-1$ crossings which are resolved
shaded-to-shaded in a monocyclic state (since in any finite tree there
is one more node than edges).  So if we replace one monocyclic state
by another then there are as many crossings whose resolution changes
from unshaded-to-unshaded to shaded-to-shaded as crossings whose
resolution changes in the other direction.  Hence there are in total
an even number of crossings whose resolution changes, proving
Lemma~\ref{Lem:even}.

\begin{Prop}
  Let $L$ be any link diagram.  Then $\langle L\rangle \in
  \capz\cdot\rootspace,$ the complex numbers of the form $pu$, where
  $p \in \capz$ and $u \in \rootspace.$ In particular $|\langle
  L\rangle |$ is a non-negative integer.
\label{Prop:pu}
\end{Prop}

\bf Proof. \rm Consider two monocyclic states for $L$.  By Lemma
\ref{Lem:even} these states disagree at an even number of crossings.
Since each change of state at a crossing introduces a factor of $A^2$
or $A^{-2}$, an even number of changes of state introduces a factor
which is a power of $A^4$ which evaluates to $-1$.

Now if $\langle L\rangle \ne 0$ then there is a monocyclic state whose
coefficient in the Kauffman bracket is a power of $A = e^{\pi
i/4} \in \rootspace$; say it is $u \in \rootspace.$ Then any other
monocyclic state has coefficient $\pm u$.  Adding together these state
values we obtain the value for $\langle L\rangle $ of $pu$ where $p
\in \capz$. This is Proposition~\ref{Prop:pu}.

The integer $|\langle L\rangle |$ is known classically as the
{\em determinant} of the link $L$.  See Section~\ref{sec:determinant}.  It
is an invariant of ambient isotopy:

\begin{Prop}
  $|\langle L\rangle |$ depends only on the ambient isotopy class of
  the link diagram $L$.
\label{Prop:determinant}
\end{Prop}

\bf Proof.  \rm  We need only check that
$|\langle L\rangle |$ remains invariant under each of the three types
of Reidemeister moves.  See Lemma~\ref{Lem:value}.  Reidemeister Move
Type~I multiplies $\langle L\rangle $ by $A^{\pm 1}=e^{\pm\pi i/4}$
and so leaves the absolute value of the bracket invariant.
Reidemeister Moves of type II and III leave the bracket, and therefore
its absolute value, unchanged.  Proposition~\ref{Prop:determinant}
follows.

\begin{Rem}
  In view of Proposition \ref{Prop:determinant} we will allow ourselves to
consider the expression $|\langle L\rangle |$ where $L$ is an ambient isotopy class of links or link diagrams, even though $\langle L\rangle $ is only
defined on regular isotopy classes of link diagrams.
\label{Rem:abuse}
\end{Rem}

\section{Numerator and Denominator}
\label{sec:numden_links}

\begin{figure}
\begin{center}
\leavevmode
\epsfbox{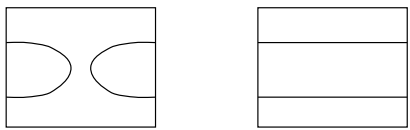}
\end{center} 
\caption{Numerator and denominator states}
\label{fig:numden}
\end{figure}

Let $T$ be a tangle.  We wish to consider what sort of links the
tangle $T$ can sit inside.  Two such links are $n(T)$ and $d(T)$,
depicted in Figure~\ref{fig:numden_links}.  If $T$ is a tangle diagram, then $n(T)$ and $d(T)$ are link diagrams.  Since the crossings of
$n(T)$ and of $d(T)$ are precisely the crossings of $T$, states for
$T$ correspond to states for $n(T)$ and to states for $d(T)$.  There
are only two homotopy types of acyclic tangle states: the {\em
numerator} states and the {\em denominator} states, as shown in
Figure~\ref{fig:numden}.  Numerator states for $T$ are then precisely
the states which correspond to monocyclic states for $n(T)$ and thus
are precisely the states which contribute to $\langle n(T)\rangle $:
$$
  \langle n(T)\rangle  = \sum_{\parbox{1.5cm}{\scriptsize\centering
  numerator\\[-5pt] states 
  for $T$}}\;\prod_{\mbox{\scriptsize vertices}} A^{\pm 1}  
$$
and similarly for $\langle d(T)\rangle $ and denominator states.  Here
by ``vertices'' we mean 4-vertices of the original tangle shadow;
ie. crossings of the tangle; for a given state at a given crossing we choose
$A$ or $A^{-1}$ according to the skein relation for the Kauffman
bracket (see \cite{state}).

\begin{Lem}
  Consider a tangle $T$.  Any numerator state for $T$ differs from any
  denominator state for $T$ at an odd number of crossings.
\label{Lem:odd}
\end{Lem}

\begin{figure}
\begin{center}
\leavevmode
\epsfbox{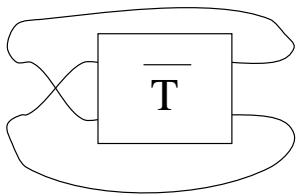}
\end{center} 
\caption{A link shadow based on the tangle shadow $\overline{T}$}
\label{fig:odd1}
\end{figure}

\bf Proof.  \rm  Consider the link shadow $L$ given in
Figure~\ref{fig:odd1}, where $\overline{T}$ is the shadow of $T$.  

\begin{figure}
\begin{center}
\leavevmode
\epsfbox{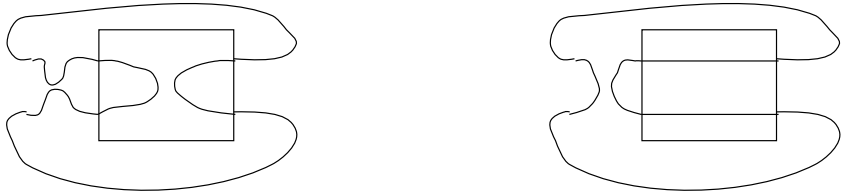}
\end{center} 
\caption{Two monocyclic link states}
\label{fig:odd2}
\end{figure}

For any numerator (denominator respectively) state for $T$ we have
that the state of Figure~\ref{fig:odd2} on the left (on the right
respectively) is monocyclic.  (The portion of the state inside the box
is shown up to homotopy.)  By Lemma~\ref{Lem:even} these states for L
differ at an even number of crossings.  Since they differ at one
crossing outside the box they must differ at an odd number of
crossings inside the box.  Hence any numerator state for T differs
from any denominator state for T at an odd number of crossings.  This
is Lemma~\ref{Lem:odd}.

\begin{Prop}
  If $\langle n(T)\rangle = pu$ and $\langle d(T)\rangle = qv$, with
  $p,q \in \capz$ and $u,v \in \rootspace$, then at least one of
  the following four equalities holds:
$$p=0,\ q=0,\ \frac{u}{v}=i,\ \frac{u}{v}=-i.$$
\label{Prop:four}
\end{Prop}

\bf Proof.  \rm If $\langle n(T)\rangle = 0$ then $pu=0$ implies that
$p=0$; if $\langle d(T)\rangle = 0$ then $qv=0$ implies that $q=0$.
In either case we are done.  If on the other hand both $\langle
n(T)\rangle $ and $\langle d(T)\rangle $ are non-zero then there is a
numerator state of coefficient, say, $u'$, and a
denominator state with coefficient, say, $v'$, where $u',v'\in\rootspace$.  By Lemma~\ref{Lem:odd} these
states differ at an odd number of crossings.  Since each crossing
change contributes a factor of $A^2 = i$ or $A^{-2} = -i$, an odd
number of crossing changes contributes a factor of $\pm i$, and we
have $u'/v' = \pm i$.  Following the proof of
Proposition~\ref{Prop:pu}, $\langle n(T)\rangle $ is of form $p'u'$
for $p'\in \capz$.  We can assume $p \ne 0$.  Then $pu=p'u'$ implies
that $p'/p= u/u'\in\capq\cap\rootspace=\{\pm 1\}$, or $u=\pm u'$.
Similarly $v= \pm v'$.  Thus $u/v = \pm u'/v' = \pm(\pm i) = \pm i$.
Proposition~\ref{Prop:four} follows.

\section{Formal Fractions}
\label{sec:range}

The invariant of tangles $f$ which we will define in the next section
takes values that are {\em formal fractions}, that is, formal quotients of two
integers.  Such fractions need not be reduced and can have zero as
their denominators.  We give a precise definition of such fractions
below.

We define an equivalence relation $\sim$ on the collection $\capz
\times \capz$ of ordered pairs of integers as follows:

\ \ \ $(p,q)\sim(p',q')$ if $(p,q)=(p',q')$ or $(p,q)=(-p',-q').$

Let us denote the set of equivalence classes as $\capa$.  We will use
the notation $[p,q]$ or $p/q$ for the equivalence class
containing the ordered pair $(p,q).$ In regards to the parallel with
rational numbers, we may think of elements of $\capa$ as fractions
which are not allowed to be reduced except by the factor $-1$.  Thus $2/4
\not= 1/2$, but $1/2 = -1/-2$ and
$3/0 = -3/0$.

We can make $\capa$ a monoid by the rule
$$\frac{p}{q} + \frac{r}{s} = \frac{ps+qr}{qs}$$
as we do with
fractions except that for the common denominator we always use the
product.

A straightforward calculation shows that this monoid structure
on $\capa$ is associative even when some of the denominators are zero.

Given the element $[p,q] \in \capa$, we will frequently speak of
properties of $p$ and $q$ for properties of the type that, if held by the pair
$(a,b)$, are also held by the pair $(-a,-b).$ Thus it makes sense to
say that, eg. ``$q$ is odd,'' ``the gcd of $p$ and $q$ is $d$,'' and
``$p$ is zero.''

\section{Definition of the Invariant}
\label{sec:definition}

In this section we will define a map  $f$ from the space $\capt$ of tangles
to the space $\capa$ of formal fractions.

\begin{Not} For a tangle diagram $T$ define the set
$$\Psi(T)=\{(u\langle n(T)\rangle ,ui\langle d(T)\rangle )|u\in\rootspace\}\subset \capc^2.$$
\end{Not}

\begin{Lem}
  The set $\Psi(T)$ depends only on the ambient
  isotopy class of the tangle diagram $T$.
\label{Lem:7ambient}
\end{Lem}

\bf Proof.  \rm We must check that $\Psi(T)$ is invariant under each of the
Reidemeister moves, for which we refer to Lemma~\ref{Lem:value}.  A
Reidemeister move on the tangle $T$ corresponds to a Reidemeister move
on $n(T)$ and to a Reidemeister move on $d(T)$.  If this move is of
type II or III the bracket evaluated at these links is unaffected;
thus $\Psi(T)$ remains invariant as claimed.

If this move is of type I, then $\langle n(T)\rangle$ and $\langle
d(T)\rangle$ are both multiplied by the same element $A^{\pm 1}\in
\rootspace$.  Thus the elements of $\Psi(T)$ are permuted and the set
itself is left unchanged.  This establishes Lemma~\ref{Lem:7ambient}.

In view of Lemma \ref{Lem:7ambient} we will refer to $\Psi(T)$ for
tangles $T$ with no particular choice of diagram.

\begin{Tad}
  Let $T\in\capt$ be a tangle. Then $\Psi(T)\cap{\capz^2}$ is a set of
  the form $\{(p,q),(-p,-q)\}$. Thus we may regard
  $\Psi(T)\cap{\capz^2}$ as an element $[p,q]$ of the set $\capa$ of
  formal fractions defined in Section~\ref{sec:range}.  This defines a
  function $f:\capt\rightarrow\capa$.
\label{Tad}
\end{Tad}

\bf Proof. \rm
Let $T$ be a fixed tangle.  We must prove three items:

\begin{enumerate}
\item $f(T)=\Psi(T)\cap{\capz^2}$ is non-empty.

\item If both $(p,q)$ and $(p',q')$ are in $f(T)=\Psi(T)\cap{\capz^2}$
then $[p,q] = [p',q'].$

\item If $[p,q] = [p',q']$ and $(p,q)\in f(T)$
 then $(p',q')\in f(T)$.
\end{enumerate}
 
\begin{description}
\item[Proof of 1]  Let $\hat{T}$ be a diagram of $T$. From Proposition
~\ref{Prop:pu} let $\langle n(\hat{T})\rangle =pu$ and $\langle
 d(\hat{T})\rangle =qv$.  From Proposition~\ref{Prop:four} we have
 $p=0$, $q=0$, $u/v =i$ or $u/v=-i.$

\begin{itemize}
\item Case $p=0$.  Set $u'=-i\overline v.$ Then $u'\langle
n(\hat{T})\rangle =u'\cdot 0=0\in\capz$ and $u'i\langle d(T)\rangle
=-i\overline{v}iqv=q\in\capz.$ Thus $(0,q)\in\Psi(T)\cap{\capz^2}$.

\item Case $q=0$.  Similar.

\item Case $u/v = i$.  Let $u' = \overline{u}$.  Then
$u'\langle n(\hat{T})\rangle =\overline u pu=p\in\capz$ and
$u'i\langle d(\hat{T})\rangle =\overline u iqv=\overline u uq=q\in\capz$.  Thus
$(p,q)\in\Psi(T)\cap{\capz^=2}$.

\item Case $u/v = -i$.  Similar.
\end{itemize}

\item[Proof of 2]  Suppose that $\hat{T}$ is a diagram of $T$ and
both $(p,q)$ and $(p',q')$ are elements of
$f(T)=\Psi(T)\cap{\capz^2}=\Psi(\hat{T})\cap{\capz^2}.$ Then there
exist elements $u$ and $u'$ of $\rootspace$ such that $p=u\langle
n(\hat{T})\rangle $, $q=ui\langle d(\hat{T})\rangle $, $p'=u'\langle
n(\hat{T})\rangle $, and $q'=u'i\langle d(\hat{T})\rangle $.

\begin{itemize}
\item Case $\langle n(\hat{T})\rangle =0$.  Then $p=p'=0$.  If $q=0$
then $\langle d(\hat{T})\rangle =0,$ which implies that $q'=0$.  So
$[p,q]=[p',q']=[0,0]$ and we are done.  If $q\ne 0$, then $\langle
d(\hat{T})\rangle \ne 0.$ We then have $q^\prime/q = u^\prime i/ui
\in \capq\cap\rootspace=\{\pm 1\}$, so that $q=\pm q'.$ Hence $[p',q']=[0,\pm
q]=[p,q],$ and we are done.

\item Case $\langle n(\hat{T})\rangle \ne 0$.\\
Then $u'\langle n(\hat{T})\rangle \ne 0.$  We have 
$u/u'=u\langle n(\hat{T})\rangle /u'\langle n(\hat{T})\rangle=p/p'\in\rootspace\cap\capq = \{\pm 1\}.$ If $u=u'$ then $p=p'$ and
$q=q'$.  Thus $[p,q]=[p',q']$.  If $u=-u'$ then $p=-p'$ and $q=-q'$.
Again $[p,q]=[p',q']$ and we are done.
\end{itemize}
 
\item[Proof of 3]  If $(p',q') \ne (p,q)$ then $(p',q')=(-p,-q);$
  replace $u$ with $-u\in\rootspace$ as the coefficient of $\langle
  n(T)\rangle$ in $\Psi(T)$.
\end{description}
This completes the proof of the first claim in
Theorem-and-Definition~\ref{Tad}; the others follow.

\begin{Prop}
If $f(T)=p/q$ then $|p|$ and $|q|$ are the determinants $|\langle n(T)\rangle
|$ and $|\langle d(T)\rangle |$, respectively.
\label{Prop:absvalue}
\end{Prop}

\bf Proof. \rm We have that $(p,q)\in\Psi(T)$ so that for some $u$ in
$\rootspace$, $p=u\langle n(T)\rangle$ and $q=ui\langle d(T)\rangle$.
So $|p|=|u\langle n(T)\rangle|=|\langle n(T)\rangle|$ and
$|q|=|ui\langle d(T)\rangle|=|\langle d(T)\rangle|$ as claimed.

\section{Properties of the Invariant}
\label{sec:properties}

\begin{Prop}
  Additivity.  If $T$ and $T'$ are tangles then $f(T+T') = f(T) +
  f(T'),$ where the addition on the left-hand side of this equation is
  defined near the end of Section~\ref{sec:domain} and the addition on
  the right-hand side is defined in Section~\ref{sec:range}.
\label{Prop:8additive}
\end{Prop}

\bf Proof. \rm  Let $T$ and $T'$ be tangle diagrams.  We have:

{\large\openup 2\jot
\begin{equation}
\begin{split}
  \left\langle n(T + T^\prime) \right\rangle  & = \hspace*{-1em}
  \sum_{\parbox{2.4cm}{\scriptsize\centering
  type \hspace*{2pt}\raisebox{3pt}{\beginpicture
        \setcoordinatesystem units <0.65pt,0.4pt>
        \putrule from -10 -10  to -10 10
        \putrule from -10 10  to 10 10
        \putrule from 10 10 to 10 -10
        \putrule from 10 -10 to -10 -10
        \ellipticalarc axes ratio 2:1 180 degrees from 10 5.6 center at 10 0
        \ellipticalarc axes ratio 2:1 -180 degrees from -10 5.6 center at -10 0
        \endpicture}\\
        acyclic states\\ for $T+T'$}}
  \!\!\!\prod_{\mbox{\scriptsize vertices}}  A^{\pm 1}\\
  & = \sum_{\parbox{2.4cm}{\scriptsize\centering type
        \hspace*{2pt}\raisebox{3pt}{\beginpicture
        \setcoordinatesystem units <0.5pt,0.4pt>
        \putrule from -10 -10  to -10 10
        \putrule from -10 10  to 10 10
        \putrule from 10 10 to 10 -10
        \putrule from 10 -10 to -10 -10
        \putrule from -10 10 to -30 10 
        \putrule from -10 -10 to -30 -10
        \putrule from -30 -10 to -30 10
        \putrule from  10 5.6 to -10 5.6
        \putrule from  10 -5.6 to -10 -5.6
        \ellipticalarc axes ratio 1.8:1 180 degrees from -10 5.6 center at -10 0
        \ellipticalarc axes ratio 1.8:1 -180 degrees from -30 5.6 center at -30 0
        \endpicture}\\
        acyclic states}}
        \!\!\!\prod_{\mbox{\scriptsize vertices}}  A^{\pm 1}
  \;\;+ \sum_{\parbox{2.4cm}{\scriptsize\centering type
        \hspace*{2pt}\raisebox{3pt}{\beginpicture
        \setcoordinatesystem units <0.5pt,0.4pt>
        \putrule from -10 -10  to -10 10
        \putrule from -10 10  to 10 10
        \putrule from 10 10 to 10 -10
        \putrule from 10 -10 to -10 -10
        \putrule from -10 10 to -30 10 
        \putrule from -10 -10 to -30 -10
        \putrule from -30 -10 to -30 10
        \putrule from -10 5.6 to -30 5.6
        \putrule from -10 -5.6 to -30 -5.6
        \ellipticalarc axes ratio 1.8:1 180 degrees from 10 5.6 center at 10 0
        \ellipticalarc axes ratio 1.8:1 -180 degrees from -10 5.6 center at -10 0
        \endpicture}\\
        acyclic states}}
        \!\!\!\prod_{\mbox{\scriptsize vertices}}  A^{\pm 1}\\
    & = \sum_{\parbox{2.4cm}{\scriptsize\centering type
  \hspace*{2pt}\raisebox{3pt}{\beginpicture
        \setcoordinatesystem units <0.65pt,0.4pt>
        \putrule from -10 -10  to -10 10
        \putrule from -10 10  to 10 10
        \putrule from 10 10 to 10 -10
        \putrule from 10 -10 to -10 -10
        \ellipticalarc axes ratio 2:1 180 degrees from 10 5.6 center at 10 0
        \ellipticalarc axes ratio 2:1 -180 degrees from -10 5.6 center at -10 0
        \endpicture}\\
        acyclic states\\ for $T$}}
  \!\!\!\prod_{\mbox{\scriptsize vertices}}  A^{\pm 1}
  \;\;\times \sum_{\parbox{2.4cm}{\scriptsize\centering type
  \hspace*{2pt}\raisebox{3pt}{\beginpicture
        \setcoordinatesystem units <0.65pt,0.4pt>
        \putrule from -10 -10  to -10 10
        \putrule from -10 10  to 10 10
        \putrule from 10 10 to 10 -10
        \putrule from 10 -10 to -10 -10
        \putrule from 10 5.6 to -10 5.6
        \putrule from 10 -5.6 to -10 -5.6
        \endpicture}\\
        acyclic states\\ for $T'$}}
  \!\!\!\prod_{\mbox{\scriptsize vertices}}  A^{\pm 1}\\
    & + \sum_{\parbox{2.4cm}{\scriptsize\centering type
  \hspace*{2pt}\raisebox{3pt}{\beginpicture
        \setcoordinatesystem units <0.65pt,0.4pt>
        \putrule from -10 -10  to -10 10
        \putrule from -10 10  to 10 10
        \putrule from 10 10 to 10 -10
        \putrule from 10 -10 to -10 -10
        \putrule from 10 5.6 to -10 5.6
        \putrule from 10 -5.6 to -10 -5.6
        \endpicture}\\
        acyclic states\\ for $T$}}
  \!\!\!\prod_{\mbox{\scriptsize vertices}}  A^{\pm 1}
  \;\;\times \sum_{\parbox{2.4cm}{\scriptsize\centering type
  \hspace*{2pt}\raisebox{3pt}{\beginpicture
        \setcoordinatesystem units <0.65pt,0.4pt>
        \putrule from -10 -10  to -10 10
        \putrule from -10 10  to 10 10
        \putrule from 10 10 to 10 -10
        \putrule from 10 -10 to -10 -10
        \ellipticalarc axes ratio 2:1 180 degrees from 10 5.6 center at 10 0
        \ellipticalarc axes ratio 2:1 -180 degrees from -10 5.6 center at -10 0
        \endpicture}\\
        acyclic states\\ for $T'$}}
  \!\!\!\prod_{\mbox{\scriptsize vertices}}  A^{\pm 1}\\
    & = \left\langle n(T) \right\rangle \langle d(T') \rangle + \langle d(T) \rangle \langle n(T') \rangle
\end{split}
\label{eq:numerator1}
\end{equation}
\begin{equation}
\begin{split}
  \left\langle d(T + T^\prime) \right\rangle  & = \hspace*{-1em}
  \sum_{\parbox{2.4cm}{\scriptsize\centering
  type \hspace*{2pt}\raisebox{3pt}{\beginpicture
        \setcoordinatesystem units <0.65pt,0.4pt>
        \putrule from -10 -10  to -10 10
        \putrule from -10 10  to 10 10
        \putrule from 10 10 to 10 -10
        \putrule from 10 -10 to -10 -10
        \putrule from 10 5.6 to -10 5.6
        \putrule from 10 -5.6 to -10 -5.6
        \endpicture}\\
        acyclic states\\ for $T+T'$}}
  \!\!\!\prod_{\mbox{\scriptsize vertices}}  A^{\pm 1}\\
  & = \sum_{\parbox{2.4cm}{\scriptsize\centering type
        \hspace*{2pt}\raisebox{3pt}{\beginpicture
        \setcoordinatesystem units <0.5pt,0.4pt>
        \putrule from -10 -10  to -10 10
        \putrule from -10 10  to 10 10
        \putrule from 10 10 to 10 -10
        \putrule from 10 -10 to -10 -10
        \putrule from -10 10 to -30 10 
        \putrule from -10 -10 to -30 -10
        \putrule from -30 -10 to -30 10
        \putrule from  10 5.6 to -10 5.6
        \putrule from  10 -5.6 to -10 -5.6
        \putrule from -10 5.6 to -30 5.6
        \putrule from -10 -5.6 to -30 -5.6
        \endpicture}\\
        acyclic states}}
        \!\!\!\prod_{\mbox{\scriptsize vertices}}  A^{\pm 1}\\
    & = \sum_{\parbox{2.4cm}{\scriptsize\centering type
  \hspace*{2pt}\raisebox{3pt}{\beginpicture
        \setcoordinatesystem units <0.65pt,0.4pt>
        \putrule from -10 -10  to -10 10
        \putrule from -10 10  to 10 10
        \putrule from 10 10 to 10 -10
        \putrule from 10 -10 to -10 -10
        \putrule from 10 5.6 to -10 5.6
        \putrule from 10 -5.6 to -10 -5.6
        \endpicture}\\
        acyclic states\\ for $T$}}
  \!\!\!\prod_{\mbox{\scriptsize vertices}}  A^{\pm 1}
  \;\;\times \sum_{\parbox{2.4cm}{\scriptsize\centering type
  \hspace*{2pt}\raisebox{3pt}{\beginpicture
        \setcoordinatesystem units <0.65pt,0.4pt>
        \putrule from -10 -10  to -10 10
        \putrule from -10 10  to 10 10
        \putrule from 10 10 to 10 -10
        \putrule from 10 -10 to -10 -10
        \putrule from 10 5.6 to -10 5.6
        \putrule from 10 -5.6 to -10 -5.6
        \endpicture}\\
        acyclic states\\ for $T'$}}
  \!\!\!\prod_{\mbox{\scriptsize vertices}}  A^{\pm 1}\\
    & = \left\langle d(T) \right\rangle \langle d(T') \rangle 
\end{split}
\label{eq:denominator1}
\end{equation}}

In equations~(\ref{eq:numerator1}) and (\ref{eq:denominator1}) we are summing
over acyclic states of the given homotopy class.  States for a sum of two tangle diagrams
correspond to pairs of states-- one for each tangle.  The double boxes
indicate the type of this pair.  

Suppose that $T$ and $T'$ take the values $[p,q]$ and $[p',q']$,
respectively.  Then there exist $u, u'\in\rootspace$ such that
$p=u\langle n(T)\rangle $, $q=ui\langle d(T)\rangle $, $p'=u'{\langle
n(T')\rangle }$, and $q'=u'i\langle d(T')\rangle .$ Substituting in
equation (\ref{eq:numerator1}), we have

\begin{equation}
\langle n(T+T')\rangle=p{\overline u}(-iq'\overline{u'})+-iq{\overline
  u}p'(\overline{u'})=-i{\overline u}{\overline{u'}}(pq'+p'q).
\label{eq:numerator2}
\end{equation}

and substituting in equation (\ref{eq:denominator1}),

\begin{equation}
\langle d(T+T')\rangle=(-iq{\overline
  u})(-iq'{\overline{u'}})=-qq'\overline u {\overline{u'}}.
\label{eq:denominator2}
\end{equation}

Multiplying (\ref{eq:numerator2}) and (\ref{eq:denominator2}) by
$iuu'$ yields
$$iuu'\langle n(T+T')\rangle =pq'+p'q$$
and
$$(iuu')i\langle d(T+T')\rangle =(iuu')i(-qq'\overline
u\overline{u'})=qq'.$$
Thus, since $iuu'\in\rootspace$, the pair
$(pq'+p'q,qq')$ belongs to the set $\Psi(T+T')$.  The value of the
invariant is $(pq'+p'q)/qq'=p/q+p'/q'.$ Proposition
\ref{Prop:8additive} follows.

\begin{figure}
\begin{center}
  \leavevmode \epsfbox{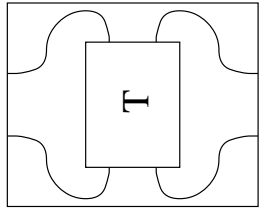}
\end{center} 
\caption{Rotation of a tangle $T$}
\label{fig:rotation_1}
\end{figure}

\begin{Not}
  If $T$ is a tangle, define $T^*$ as the tangle obtained by rotating
  $T$ by the angle $\pi/2$ and redirecting the endpoints (see
  Figure~\ref{fig:rotation_1}).
\label{not:rotation}
\end{Not}

\begin{Prop}
  If $T$ takes the value $[p,q]$ then $T^*$ takes the value $[-q,p]$.
\label{Prop:8rotation}
\end{Prop}

\bf Proof. \rm  First note that given a diagram for $T$ Figure~\ref{fig:rotation_1} yields a diagram for $T^*$.  Let us fix a diagram $D$ for $T$ and let $D^*$ be the corresponding diagram for $T^*$.

\begin{figure}
\begin{center}
  \leavevmode \epsfbox{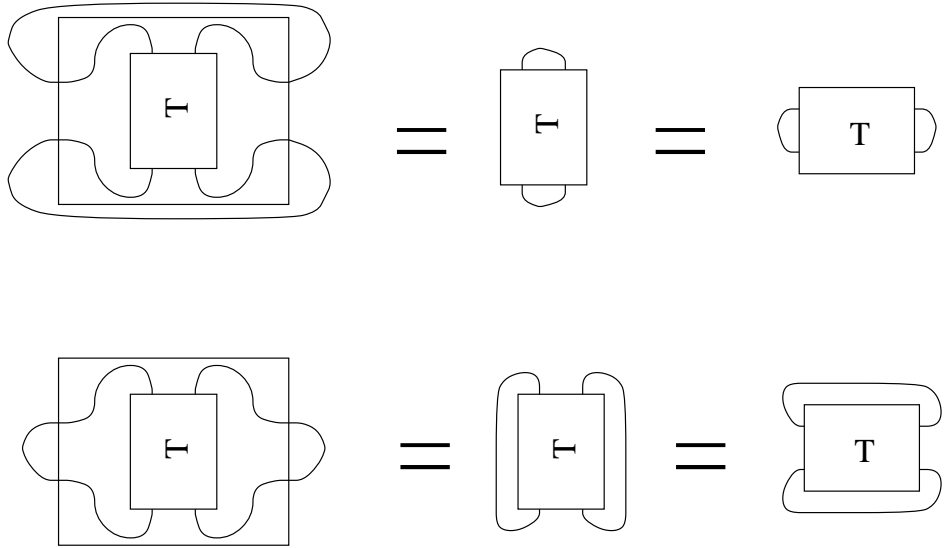}
\end{center} 
\caption{Isotopy of numerator and denominator closures}
\label{fig:rotation_2}
\end{figure}

The first line of Figure \ref{fig:rotation_2} shows that $n(D^*)$ is
     regular isotopic to $d(D)$ and the second line shows that
     $d(D^*)$ is regular isotopic to $n(D)$.  Since $T$ takes the
     value $[p,q]$, the pair $(p,q)\in\Psi(D),$ ie. there exists a
     $u\in\rootspace$ such that $p=u\langle n(D)\rangle $ and
     $q=iu\langle d(D)\rangle ,$ or $\langle n(D)\rangle = p\overline
     u$ and $\langle d(D)\rangle =q(-i\overline u)$.

Set $u^*=iu\in \rootspace$.

Then

$$u^*\langle n(D^*)\rangle =u^*\langle d(D)\rangle =iuq(-i\overline u)=q,$$ and
$$u^*i\langle d(D^*)\rangle =u^*i\langle n(D)\rangle =iuip\overline u =-p.$$

Thus $u^*\in\rootspace$ implies that $(q,-p)\in\Psi(D^*)$; the value of
the invariant is $[q,-p]=[-q,p]$.  This completes the proof of
Proposition~\ref{Prop:8rotation}.  

\begin{Not}
  In view of Proposition \ref{Prop:8rotation} we will introduce an
  operation labelled $*$ on $\capa$, defined by $[p,q]^*=[-q,p].$
  Then Proposition~\ref{Prop:8rotation} asserts that $f(T^*)=f(T)^*$
  for all tangles $T\in\capt$.
\label{Not:star}
\end{Not}

\begin{Prop}
  Reflecting a tangle $(J,T)$ about any of the three planes $x=1/2$,
  $y= 1/2$, or $z=1/2$ yields another tangle $(J,\overline T)$ whose
  invariant is $-p/q$, where $p/q$ is the value of the invariant
  on $T.$
\label{Prop:8reflect}
\end{Prop}

\bf Proof. \rm Let $D$ be a diagram for the tangle $T$, and $\overline
D$ be the reflection of $D$ through the given plane.  If $S$ is a
state for $D$ then $S$ may also be reflected through the given plane
(in the case of the plane $z=1/2$ the state $S$ is left unchanged)
yielding a state $\overline S$ for $\overline D$.  Note that
$\overline S$ is a numerator (denominator) state precisely if $S$ is.
Therefore the terms in the expansion of the Kauffman bracket of
$n(\overline D)$ correspond precisely to those of $n(D)$ and similarly
for $d(\overline D)$ and $d(D)$.  Since crossings for $D$ also
correspond to crossings for $\overline D$ the factors of each term of
the expansion for $\langle n(\overline D)\rangle$ correspond precisely
to those for $\langle n(D)\rangle$ and similarly for $\langle
d(\overline D) \rangle$ and $\langle d(D)\rangle$.  But if $S$ is resolved as an $A$-state at a
crossing of $D$ then $\overline S$ is resolved as an $A^{-1}$-state at
the corresponding crossing of $\overline D$ and similarly with the
roles of $A$ and $A^{-1}$ reversed.  By $A^{-1}=\overline A$ (here
$\overline A$ is the complex conjugate of $A$) and the fact that the
Kauffman bracket is a polynomial in $A$ with real (in fact integer)
coefficients we have that $\langle n(\overline D)\rangle
=\overline{\langle n(D)\rangle }$ and $\langle d(\overline D)\rangle
=\overline{\langle d(D)\rangle }$.

As before, write $\langle n(D)\rangle  = p\overline u$ and $\langle d(D)\rangle =q(-i\overline u)$ for some $u\in\rootspace.$

We have

$$\langle n(\overline D)\rangle =\overline{\langle n(D)\rangle
}=\overline{p\overline u}=pu$$
$$\langle d(\overline D)\rangle =\overline{\langle d(D)\rangle }=\overline{-iq\overline u}=qui.$$

Set $u' = -\overline u.$  Then
$$u'\langle n(\overline D)\rangle =-{\overline u}pu=-p$$ and
$$u'i\langle d(\overline D)\rangle =-{\overline u}i(qui)=q.$$

Thus $u'\in\rootspace$ implies that $(-p, q)\in\Psi(\overline D)$.
The value of the invariant on the tangle $\overline T$ is $[-p,q]$.
This is Proposition~\ref{Prop:8reflect}.

\begin{Prop}
If $\langle L\rangle = pu$, with $p \in \capz$ and $u \in \rootspace$, then $p$ (or the determinant $|p|$) is even if
and only if $L$ has more than one component.
\end{Prop}

\bf Proof. \rm It will be shown in Section \ref{sec:determinant} that
$p=\pm V_L (-1)$, where $V$ is the Jones polynomial.  The result now
follows from \cite[Thm. 15, p. 107]{jones} noting that $\pm V_L (-1)
\equiv V_L (1)$ mod $2$).

From Proposition~\ref{Prop:absvalue} we now have:

\begin{Prop}
For a 4-tangle $T$ with invariant
$[p,q]$ the parities of $p$ and $q$ are related to the homotopy class
of $T$ according to Figure~\ref{fig:htpy_table}.
\label{Prop:pqparity}
\end{Prop}

\begin{figure}
\font\thinlinefont=cmr5
\begingroup\makeatletter\ifx\SetFigFont\undefined%
\gdef\SetFigFont#1#2#3#4#5{%
  \reset@font\fontsize{#1}{#2pt}%
  \fontfamily{#3}\fontseries{#4}\fontshape{#5}%
  \selectfont}%
\fi\endgroup%
\mbox{\beginpicture
\setcoordinatesystem units <1.00000cm,1.00000cm>
\unitlength=1.00000cm
\linethickness=1pt
\setplotsymbol ({\makebox(0,0)[l]{\tencirc\symbol{'160}}})
\setshadesymbol ({\thinlinefont .})
\setlinear
%
%
\linethickness= 0.500pt
\setplotsymbol ({\thinlinefont .})
\putrectangle corners at  2.826 25.622 and 10.954 20.288
%
%
\linethickness= 0.500pt
\setplotsymbol ({\thinlinefont .})
\putrectangle corners at  2.826 25.622 and  7.144 22.955
%
%
\linethickness= 0.500pt
\setplotsymbol ({\thinlinefont .})
\putrectangle corners at  7.144 22.955 and 10.954 20.288
%
%
\linethickness= 0.500pt
\setplotsymbol ({\thinlinefont .})
\putrectangle corners at  7.779 25.114 and  9.557 23.844
%
%
\linethickness= 0.500pt
\setplotsymbol ({\thinlinefont .})
\putrectangle corners at  3.588 22.447 and  5.366 21.177
%
%
\linethickness= 0.500pt
\setplotsymbol ({\thinlinefont .})
\putrule from  3.588 22.066 to  5.366 22.066
%
%
\linethickness= 0.500pt
\setplotsymbol ({\thinlinefont .})
\putrule from  3.588 21.558 to  5.366 21.558
%
%
\linethickness= 0.500pt
\setplotsymbol ({\thinlinefont .})
\putrectangle corners at  7.779 22.447 and  9.557 21.177
%
%
\linethickness= 0.500pt
\setplotsymbol ({\thinlinefont .})
\plot  7.779 21.558  9.557 22.066 /
%
%
\linethickness= 0.500pt
\setplotsymbol ({\thinlinefont .})
\plot  7.779 22.066  8.414 21.939 /
%
%
\linethickness= 0.500pt
\setplotsymbol ({\thinlinefont .})
\plot  8.922 21.685  9.557 21.558 /
\linethickness= 0.500pt
\setplotsymbol ({\thinlinefont .})
%
%
%
\plot    7.779 24.733  7.842 24.733
         7.937 24.733
         8.009 24.733
         8.096 24.733
         8.188 24.729
         8.271 24.717
         8.346 24.698
         8.414 24.670
         8.509 24.590
         8.541 24.479
         8.525 24.368
         8.477 24.289
         8.382 24.241
         8.311 24.229
         8.223 24.225
         8.132 24.225
         8.049 24.225
         7.973 24.225
         7.906 24.225
         /
\plot  7.906 24.225  7.779 24.225 /
\linethickness= 0.500pt
\setplotsymbol ({\thinlinefont .})
%
%
%
\plot    9.557 24.733  9.493 24.733
         9.398 24.733
         9.327 24.733
         9.239 24.733
         9.148 24.729
         9.065 24.717
         8.989 24.698
         8.922 24.670
         8.826 24.590
         8.795 24.479
         8.811 24.368
         8.858 24.289
         8.954 24.241
         9.025 24.229
         9.112 24.225
         9.204 24.225
         9.287 24.225
         9.362 24.225
         9.430 24.225
         /
\plot  9.430 24.225  9.557 24.225 /
%
%
\put{\SetFigFont{20}{24.0}{\rmdefault}{\mddefault}{\updefault}$p$ even} [lB] at  3.969 26.003
%
%
\put{\SetFigFont{20}{24.0}{\rmdefault}{\mddefault}{\updefault}$p$ odd} [lB] at  8.0 26.003
%
%
\put{\SetFigFont{20}{24.0}{\rmdefault}{\mddefault}{\updefault}$q$ even} [lB] at  0.286 24.098
%
%
\put{\SetFigFont{20}{24.0}{\rmdefault}{\mddefault}{\updefault}$q$ odd} [lB] at  0.413 21.685
%
%
\put{\SetFigFont{12}{14.4}{\rmdefault}{\mddefault}{\updefault}not loop-free} [lB] at  3.969 24.225
%
%
\put{\SetFigFont{12}{14.4}{\rmdefault}{\mddefault}{\updefault}loop-free} [lB] at  8.668 23.209
%
%
\put{\SetFigFont{12}{14.4}{\rmdefault}{\mddefault}{\updefault}loop-free} [lB] at  4.477 20.542
%
%
\put{\SetFigFont{12}{14.4}{\rmdefault}{\mddefault}{\updefault}loop-free} [lB] at  8.668 20.542
\linethickness=0pt
\putrectangle corners at  0.159 26.452 and 10.979 20.263
\endpicture}
\caption{Homotopy class and parity of $p$ and $q$}
\label{fig:htpy_table}
\end{figure}

Here the diagrams in the three lower and right-hand entries describe
the homotopy class (Definition \ref{Def:acyclic}) of the tangle.  The
reader may wish to corroborate Figure~\ref{fig:htpy_table} against
Figure~\ref{fig:various} and the values of the invariant given there.

\section{The Embeddability Condition}
\label{sec:embed}

We say that the tangle $T$ can be {\em embedded in} or {\em sits
inside} the link $L$ if there is a representative of $L$ whose
intersection with the unit cube $J$ is the 1-manifold $T$.

The main result of this paper is the following:

\begin{Thm}  Suppose that the tangle $T$ can be embedded in the link
$L$.\leavevmode\\
{\bf Version 1.}\quad
If $f(T)=p/q$ then $gcd(p,q)\ \Bigm| \ 
{|\langle L\rangle |}$ as integers.

\noindent
{\bf Version 2.}\quad
If $n(T)$ and $d(T)$ are the links of Figure~\ref{fig:numden_links}, then
$${gcd(|\langle n(T)\rangle |,|\langle d(T)\rangle |)}\ \Bigm| \ 
{|\langle L\rangle |}$$ as integers.
\label{Thm:9main}
\end{Thm}

Since the determinant of the unknot is $1$, we have

\begin{Cor}
  The tangle $T$ cannot sit inside the unknot unless the determinants
$|\langle n(T)\rangle |$ and $|\langle d(T)\rangle |$ of the links $n(T)$ and $d(T)$ are relatively prime.
\label{Cor:unknot}
\end{Cor}

The proof of the theorem uses the following lemma:

\begin{Lem}
  If the tangle $T$ can be embedded in the link $L$, then there is a
  second tangle $T'$ such that $L$ is ambient isotopic to the link
  $n(T+T')$.
\label{Lem:9ambient}
\end{Lem}

\bf Proof of Lemma \ref{Lem:9ambient}.  \rm The tangle $T$ is a
properly embedded 1-manifold in a topological 3-ball (actually a cube)
$J$ contained in $S^3$ with $T\cap \partial J$ consisting of 4 points.
Let $B$ denote the closure of $S^3\setminus J$, also a 3-ball.  Then
the intersection of $L$ with $B$ is another tangle, that we can call
$T'$.  Since $T'$ is not space-filling, it is clear that $T'$ can be
brought into a cube in $\bf R^3$ adjacent to the cube $J$ by an
ambient isotopy that is the identity on the boundary of $B$.  This
produces the standard picture in the form $n(T+T')$.

\bf Proof of Theorem \ref{Thm:9main}. \rm Version 1:  Let $T'$ be the tangle
given by Lemma~\ref{Lem:9ambient}.  Write $f(T') = [r,s]$.  Then
$f(T+T')=[ps+qr,qs]$ and so by Proposition~\ref{Prop:absvalue},
$|ps+qr| = |\langle n(T+T')\rangle |=|\langle L\rangle |$ and
$$gcd(p,q)\ \Bigm| \ {|ps+qr|} = {|\langle L\rangle |}$$ as claimed.

Theorem \ref{Thm:9main}, Version 2 follows from Proposition~\ref{Prop:absvalue} and Theorem \ref{Thm:9main}, Version~1.

\section{Sample Calculations}
\label{sec:calculations}

\begin{figure}
\begin{center}
\leavevmode
\epsfbox{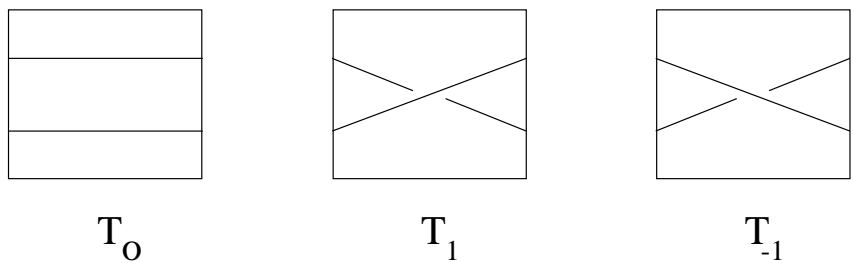}
\end{center} 
\caption{Three very simple tangles}
\label{fig:simple_tangles}
\end{figure}

For the tangle $T_0$ depicted in Figure \ref{fig:simple_tangles}, $n(T_0)$ is the unlink of two
components and $d(T_0)$ is the unknot.  Let $u=-i$.  Then
$u\langle n(T_0)\rangle =-i\cdot 0=0\in\capz$ and $ui\langle d(T_0)\rangle =\langle d(T_0)\rangle =1\in\capz$.  The
value of the invariant on $T_0$ is $0/1$.

Now consider the tangle $T_1$ of Figure \ref{fig:simple_tangles}.  We
have that $n(T_1)$ and $d(T_1)$ are both unknots and by applying the
skein relation for the Kauffman bracket as in equation (\ref{eq:Rm_1}),
$\langle n(T_1)\rangle =A$ and $\langle d(T_1)\rangle =A^{-1}$.  Set
$u= \overline{A}=e^{-\pi i/4}$.  Then $u\langle n(T_1)\rangle
=e^{-\pi i/4}e^{\pi i/4}=1\in\capz$ and $ui\langle
d(T_1)\rangle =e^{-\pi i/4}ie^{-\pi i/4}=1\in\capz$.  Thus
the value of the invariant on $T_1$ is $1/1$.

The value of the invariant on $T_{-1}$ of Figure
\ref{fig:simple_tangles} is $1/-1=-1/1$.  This value may be found by
noticing that $T_{-1}$ can be obtained from $T_1$ by either a rotation
(Proposition \ref{Prop:8rotation}) or a reflection
(Proposition~\ref{Prop:8reflect}).  Alternatively one may find this
value by solving for $f(T_{-1})$ in the formal fraction equation $1/1
+ f(T_{-1})= 0/1$ which results from applying additivity to the tangle
equation $T_1 + T_{-1}=T_0$.

\begin{figure}
\begin{center}
\leavevmode
\epsfbox{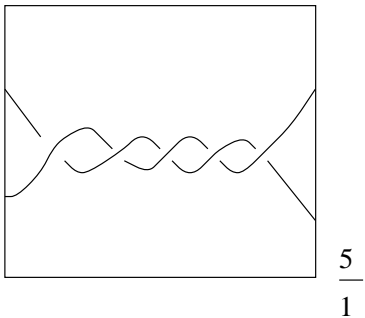}
\end{center} 
\caption{An integral tangle}
\label{fig:integral}
\end{figure}

\begin{Def}
Let $n\in\capz$; we will define the tangle $T_n$.  $T_n$ has been
defined above for $n\in\{-1,0,1\}$.  For $n>1$ define $T_n$ as the sum
of $n$ copies of $T_1;$ for $n<-1$ define $T_n$ as the sum of $|n|$
copies of $T_{-1}$.  We will call $T_n$ for $n\in\capz$ an {\em
integral} tangle as in \cite{conway}.  By the additivity of the
invariant and the associativity of the monoid structure on $\capa$,
$T_n$ is assigned the value $[n,1]=n/1$.  $T_5$ is depicted in Figure
\ref{fig:integral}.
\label{Def:integral}  
\end{Def}

\begin{Prop}
Define the function $\phi :\capz \rightarrow \capa:$ $ n \mapsto
n/1$.  Then we have that $T_{m}+T_{n}=T_{m+n}$ for all
$m,n\in\capz$ and that $\phi ^{-1} \circ f:T_n\mapsto n$ is a
well-defined group isomorphism from the submonoid of integral tangles
to the integers.
\end{Prop}

\bf Proof. \rm  The proof is entirely algebraic except for the fact that $T_n
+ T_{-n}=T_0$, by successive applications of the Reidemeister type II
move, and is left to the interested reader.

\begin{figure}
\begin{center}
\leavevmode
\epsfbox{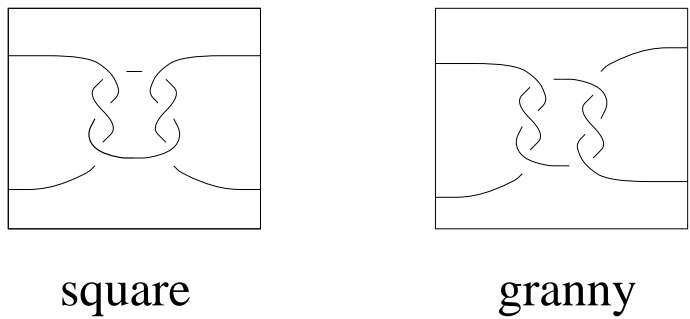}
\end{center} 
\caption{Two commonly confused tangles}
\label{fig:square_granny}
\end{figure}

\begin{Ex}
``Square knot'' tangle.

The tangle $(T_3)^*+(T_{-3})^*$, commonly known as a square knot, is
drawn on the left side of Figure~\ref{fig:square_granny}.  The value of
the invariant (see Proposition~\ref{Prop:8rotation} and Notation
 ~\ref{Not:star}) is:
  $$(\frac{3}{1})^* +{(\frac{-3}{1})^*} =\frac{-1}{3}+\frac{1}{3}=\frac{0}{9}.$$
  
  Note here that $gcd(0,9)=9$, so this tangle does not sit inside an
  unknot or for example a trefoil (determinant $=3$).
\end{Ex}

\begin{Ex}
  A ``granny knot'' tangle $G$ is depicted on the right side of
  Figure~\ref{fig:square_granny}.
$$G=(T_3)^*+(T_3)^*;$$
$$f(G)=(\frac{3}{1})^* + {(\frac{3}{1})^*}
=\frac{-1}{3}+\frac{-1}{3}=\frac{-6}{9}.$$

Here $gcd(-6,9)=3$, so $G$ possibly sits inside a
trefoil\footnote{In fact John Wood has pointed out to the author that
  it does: $n(G+T_1)=n(T_{-3})$, a right-handed trefoil.} but not, for
  instance, an unknot or a Hopf link (determinant $=2$).
\end{Ex}

\begin{figure}
\begin{center}
\leavevmode
\epsfbox{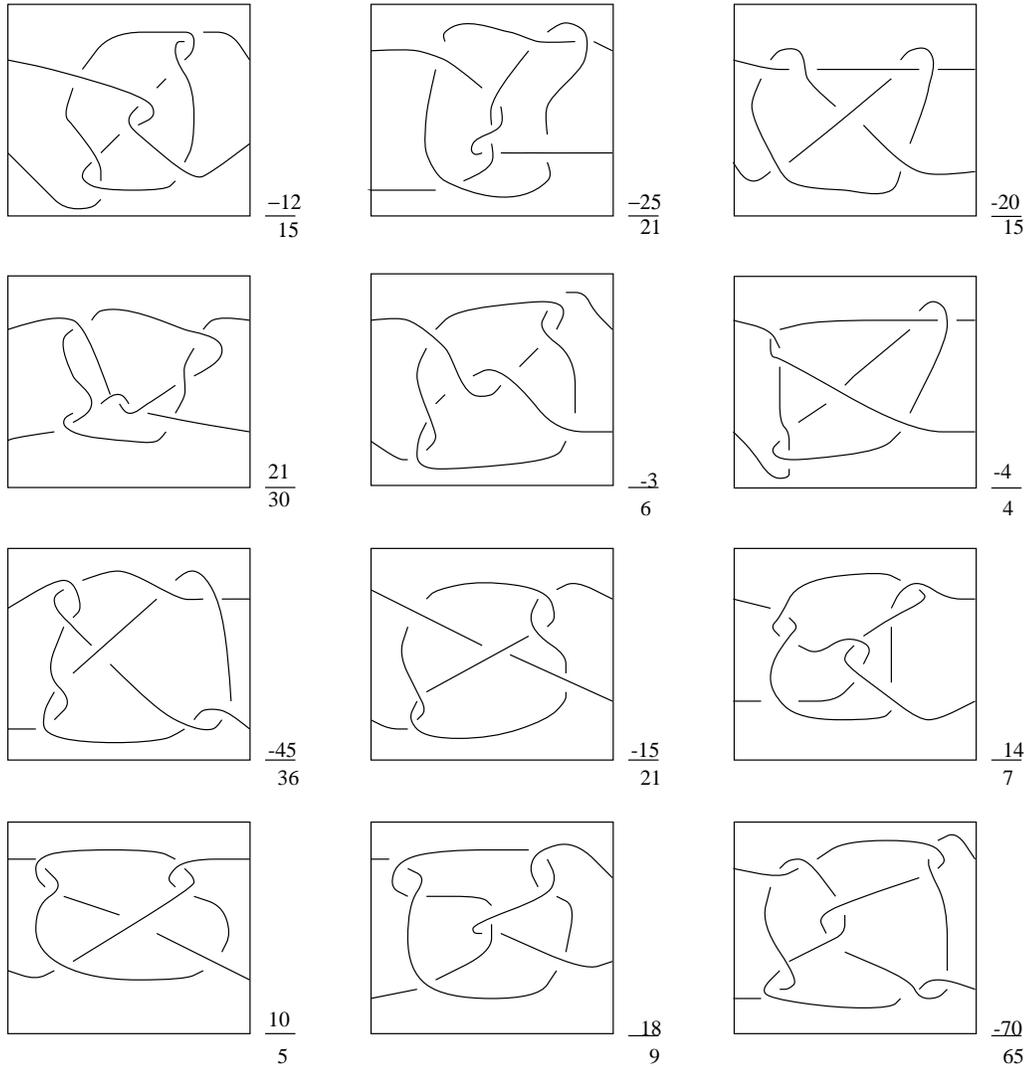}
\end{center} 
\caption{Some tangles and their associated values.}
\label{fig:various}
\end{figure}

Figure \ref{fig:various} gives the value of the invariant for various other tangles.

\section{Determinant}
\label{sec:determinant}

The non-negative integer $|\langle L\rangle |$ has many other
formulations.  Firstly, it is, up to sign, the value at $t=-1$ of the
Jones polynomial $V_L(t)$: This is because to compute the Jones
polynomial from the Kauffman bracket we replace $A^{-4}$ with $t$ and
multiply by a power of $-A^{-3}$ which has absolute value 1 for our
choice of $A$.  Since $A^{-4} = -1$, setting $t=-1$ in the Jones
polynomial will give, in absolute value, the same quantity as the
Kauffman bracket.  But the Jones polynomial has integer coefficients,
so both $V_L(-1)$ and $|\langle L\rangle |$ are integers.  Therefore
they agree up to sign.  Secondly, when $L$ is a knot $|\langle
L\rangle |$ is the value at $-1$ of the Alexander polynomial (see
\cite[Cor. 13, p. 107]{jones}).  Finally, when $L$ is a knot $|\langle
L\rangle |$ is the order of the first homology of the two-fold cover
of $S^3$ branched over $L$ (see \cite{rolfsen}).

The integer $|\langle L\rangle |$ is known classically as the {\em
determinant} of the link $L$.

\bf Remark.  \rm Probably the easiest way to
systematically calculate the value of the determinant on a large
number of links with few crossings (at most ten, say) is by
calculating the value of the invariant on a collection of related
tangles.  For instance we were able to use state summations for the
Kauffman bracket to evaluate the invariant on each of the twelve
tangles in Figure~\ref{fig:various}.  This gives the value of the
determinant for $2\cdot 12=24$ links representing at least 18 (the
number of positive integers which appear up to sign in the formal
fractions) link types, each with no more than ten crossings.

For determinants of some simple individual links, the reader is
offered two strategies:

1) For algebraic links-- links arising as $n(T)$ for some
algebraic tangle $T$ (see Section~\ref{sec:realiz}) such as the square or
granny tangle-- calculate the determinant as the absolute
value of the numerator of the invariant associated to the tangle,
which can be computed using the additivity, rotation and reflection
formulae (Section~\ref{sec:properties}).

2) Consult the tables of knot- and link-diagrams in \cite{rolfsen} (where it is
explained how to calculate the determinant from the data accompanying
each link-diagram).

\section{Realizability}
\label{sec:realiz}

Let $\tzero \subset \capt$ be the set of algebraic tangles with
unknotted components, where the set of ``algebraic'' tangles is the
closure under addition, rotation and reflection of the set of integral
tangles as in \cite{conway}, and ``unknotted components'' means that
each component of $T$ is unknotted (for the two arc components this
means that the knot formed by adjoining an arc in $\partial J$ with
the same endpoints is an unknotted loop).

We will show
\begin{Thm}
  The function f maps $\tzero$ onto $\capa$.
\label{Thm:12onto}
\end{Thm}

This defines a large (infinite) class of tangles for which Theorem~\ref{Thm:9main} gives non-trivial information.

\bf Proof of Theorem \ref{Thm:12onto}.  \rm Let $[p,q] \in \capa$, and write
$p={2^n}dp'$, $q={2^n}dq'$, where $n\ge 0$, $d$ is odd, and $p'$ and
$q'$ are relatively prime.  Note that if $q'$ is not odd, then $p'$ is
odd.  If we can show that the value $[q,p]$ is achieved on $\tzero$
then by rotating and reflecting (section~\ref{sec:properties}), we
arrive at a tangle, still in $\tzero$, which is mapped to $[p,q]$.
Thus we may assume that $q'$ is odd.

We have

\begin{equation}
\frac{p}{q} = \frac{{2^n}dp'}{{2^n}dq'} = \frac{0}{2^n}+\frac{0}{d}
+\frac{p'}{q'}.
\label{eq:factorize}
\end{equation}

We will explicitly find tangles in $\tzero$ whose values are each of
the three terms in the rightmost expression of (\ref{eq:factorize}) and
then show that their sum also lies in $\tzero$.

First, we describe the tangle on which $f$ takes value $[p',q']:$

Expand the rational number $p'/ q'$ as a finite continued fraction.
Recall that a finite continued fraction is obtained by starting with
the integer $0$ and alternately applying for a finite number of
iterations the operations of a) addition of a non-zero integer $k$,
(not necessarily the same integer $k$ for different iterations) and b)
taking the reciprocal $r \mapsto 1/r$.  Corresponding to these
algebraic operations are the topological operations of a$^\prime$) the
addition of an integral tangle $T_k$, $k\ne 0$, and b$^\prime$) a
rotation and a reflection (see Section~\ref{sec:properties}).

We notice that i) the invariant on tangle $T_0$ is $0$ and ii) an
operation of type a$^\prime$ (b$^\prime$, respectively) on a tangle
$T$ has the effect of performing an operation of type a (b,
respectively) on the invariant of $T$.  Thus we have only to expand
the rational number $p'/q'$ as a continued fraction and then
replace the operations a) and b) with the operations a$^\prime$) and
b$^\prime$) (starting from the tangle $T_0$) to arrive at tangle whose
invariant is $p'/q'$.  A tangle constructed in such a way is
called a {\em rational} tangle (see \cite{conway}).  By virtue of this
mode of construction, rational tangles are algebraic and have
unknotted strands (there are no loop components); ie. lie in $\tzero$.

Next, we describe below the tangle on which $f$ takes the value
$[0,d]$.  Let $d'$ be the integer (not the formal fraction)
$(d-1)/2.$ Consider the following calculation:

$$\frac{1}{d} + \frac{1}{d} = \frac{2d}{d^2}$$

$$(\frac{2d}{d^2})^* = \frac{-d^2}{2d}$$

$$\frac{-d^2}{2d} + \frac{d'}{1} = \frac{-d}{2d}$$

$$(\frac{-d}{2d})^*=\frac{2d}{d}$$

$$\frac{2d}{d}+\frac{-2}{1}=\frac{0}{d}.$$

Corresponding to these algebraic operations are the following
topological ones: We start with a variant of the granny ($d=3$) tangle
and alternately add integral tangles and rotate.  These operations
preserve the properties that the tangle is algebraic and that both
strands (there are no loop components since $d$ is odd) are unknotted.
By the properties of the invariant described in Section
\ref{sec:properties}, the resulting tangle $T$ satisfies
$f(T)=0/d.$  

Here is an explicit formula for $T$ of the above paragraph:

$$(((T_{-d})^*+(T_{-d})^*)^*+T_{d'})^*+T_{-2}.$$

\begin{figure}
\begin{center}
\leavevmode
\epsfbox{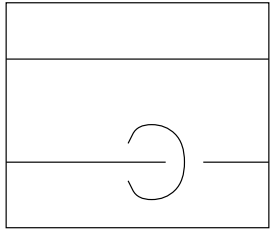}
\end{center} 
\caption{The tangle $S$ used to show realizability of formal fractions}
\label{fig:zeroovertwo}
\end{figure}

Finally, define the tangle $S$ as $({T_2^*}+{T_0^*})^*$, depicted in
Figure~\ref{fig:zeroovertwo}.  Then
$$f(S)=({(\frac{2}{1})^*}+{(\frac{0}{1})^*})^*=(\frac{-1}{2}+\frac{1}{0})^*=(\frac{2}{0})^*=\frac{0}{2}.$$
The sum of $n$ copies of $S$ lies in $\tzero$ and has invariant
$$\overbrace{\frac{0}{2} + \dots + \frac{0}{2}}^{\text{$n$ times}} = \frac{0}{2^n}.$$

It remains to be seen that the sum of these three tangles, all in
$\tzero$, also lies in $\tzero$.

Since adding multiple copies of the tangle $S$ preserves both the
properties of being algebraic and having unknotted components, we will
be done proving Theorem~\ref{Thm:12onto} if we can show that the sum of the two tangles corresponding
to the last two terms in (\ref{eq:factorize}) is in $\tzero$.  For this we
need the following lemma:

\begin{Lem}
i) Let $T$ be a tangle mapped to some element $[p,q]$ of $\capa$.  Then
$q$ is odd precisely if $T$ is loop-free and the homotopy class of $T$
is in the set 

\begin{equation}
H=\{\;\;
{\beginpicture
        \setcoordinatesystem units <0.7pt,0.56pt>
        \putrule from -10 -5   to -10 15
        \putrule from -10 15  to 10 15
        \putrule from 10 15 to 10 -5  
        \putrule from 10 -5  to -10 -5 
        \putrule from  10 10.6 to -10 10.6
        \putrule from  10 -0.6 to -10 -0.6
        \endpicture}\;\;,\;\;
{\beginpicture
        \setcoordinatesystem units <0.7pt,0.56pt>
        \putrule from -10 -5   to -10 15
        \putrule from -10 15  to 10 15
        \putrule from 10 15 to 10 -5 
        \putrule from 10 -5  to -10 -5 
        \plot  10 10.6 -10 -0.6 /
        \plot  -10 10.6 -2.5 6.4 /
        \plot  10 -0.6 2.5  4.4 /
        \endpicture}\;\;
\}.
\end{equation}

ii) Consider a pair of tangles with the following properties:
algebraic, unknotted components, loop-free, and homotopy class in $H$.
Then their sum also has these properties.
\label{Lem:12onto}
\end{Lem}

\bf Proof of Lemma \ref{Lem:12onto}.  \rm
Part i):  See Proposition~\ref{Prop:pqparity}.
Part ii):  Obvious.  The sum can be seen to have the last three
properties by examining its components one at a time.

Proof of Theorem \ref{Thm:12onto} (resumed).   We have two algebraic
tangles with unknotted components, one of which is mapped to $[0,d]$ and
the other of which is mapped to $[p',q']$.  Note that the denominators
$d$ and $q'$ are both odd.  Therefore Lemma~\ref{Lem:12onto} part i)
applies to both of these tangles.  Lemma~\ref{Lem:12onto} part ii)
then asserts that their sum is also algebraic with unknotted
components, ie.  belongs to $\tzero$ as claimed.

This completes the proof of Theorem \ref{Thm:12onto}.

\bf Remark.  \rm In the construction of the rational tangle with
invariant $[p',q']$ above there are many ways of expanding the
rational number $p'/q'$ as a continued fraction.  Furthermore, in the
construction of the rational tangle itself we may add the integral
tangle either on the left or on the right, we may rotate through an
angle either of $\pi /2$ or of $-\pi /2$ (which may be shown to have
the same effect on the invariant) and we may reflect about any of the
planes $x=1/2$, $y=1/2$, or $z=1/2$.  Thus our procedure gives many
distinct tangle diagrams associated to the same value of the
invariant.  However Goldman and Kauffman prove in \cite{rational} that
the rational tangles represented by these diagrams are all ambient
isotopic to each other.

\section{Fox $n$-colourings}
\label{sec:colourings}

\begin{figure}
\begin{center}
\leavevmode
\epsfbox{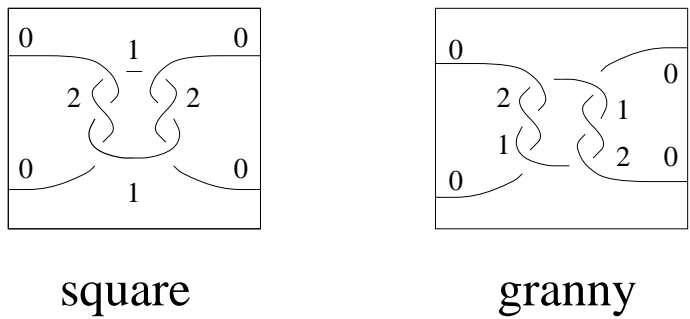}
\end{center} 
\caption{Tricolourings which show neither square nor granny tangle embeds in an unknot}
\label{fig:tricolourings}
\end{figure}

Daniel Silver has pointed out to the author that all of the tangles in
Figure~\ref{fig:various} with $gcd(p,q)$ odd and not equal to $1$, as
well as the square and granny tangles, can alternatively be shown by
Fox $n$-colourings not to be embeddable in unknots.  (For the idea,
without the terminology, of a Fox $n$-colouring see Chapter 10 of \cite{trip};
\cite{silverwms} contains an up-to-date review of this concept.)  One
need only produce a non-trivial $n$-colouring of a diagram for the
tangle in which all of the overcrossing arcs incident to one (or more)
of the four endpoints are coloured the same, say with the label
(``colour'') $a$.  Fox $n$-colourings of this type with $n=3$
(``tricolourings'') for the square and the granny tangles are depicted
in Figure~\ref{fig:tricolourings}.  If the tangle diagram could be
embedded in an unknot diagram then by colouring all the remaining arcs
(outside the tangle diagram) with the same label $a$, we would have a
non-trivial colouring of the unknot, a contradiction.  

If both the numerator and the denominator closure of a $4$-tangle $T$
can be $n$-coloured non-trivially, can the tangle be non-trivially
$n$-coloured so that all four endpoints get the same colour?  It is
known that a knot with determinant $p$ can be non-trivially
$q$-coloured, where $q$ is a prime, if and only if $q$ divides $p$.
(See \cite{livingston}, p. 45.)  If the question above were true then
we would have the following situation: If there were a prime $q$
dividing both the determinant of $n(T)$ and of $d(T)$ then the tangle
$T$ would have a non-trivial $q$-colouring in which all four endpoints
receive the same colour and thus could not sit inside an unknot.  This
would provide an alternate proof of Corollary~\ref{Cor:unknot}.

\section{Generalizations and Further Questions}
\label{sec:question}

\begin{figure}
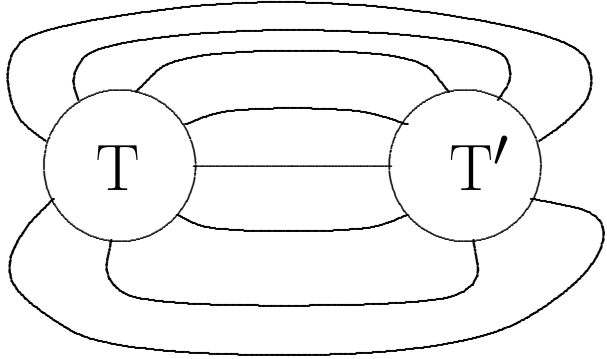

\begin{center}
\font\thinlinefont=cmr5
\begingroup\makeatletter\ifx\SetFigFont\undefined%
\gdef\SetFigFont#1#2#3#4#5{%
  \reset@font\fontsize{#1}{#2pt}%
  \fontfamily{#3}\fontseries{#4}\fontshape{#5}%
  \selectfont}%
\fi\endgroup%
\mbox{\beginpicture
\setcoordinatesystem units <1.00000cm,1.00000cm>
\unitlength=1.00000cm
\linethickness=1pt
\setplotsymbol ({\makebox(0,0)[l]{\tencirc\symbol{'160}}})
\setshadesymbol ({\thinlinefont .})
\setlinear
%
%
\linethickness= 0.500pt
\setplotsymbol ({\thinlinefont .})
\ellipticalarc axes ratio  1.008:1.008  360 degrees 
        from  7.292 24.323 center at  6.284 24.323
%
%
\linethickness= 0.500pt
\setplotsymbol ({\thinlinefont .})
\ellipticalarc axes ratio  1.008:1.008  360 degrees 
        from  2.699 24.323 center at  1.691 24.323
\linethickness= 0.500pt
\setplotsymbol ({\thinlinefont .})
%
%
%
\plot    2.673 24.323  3.988 24.323
         /
\plot  3.988 24.323  5.302 24.323 /
\linethickness= 0.500pt
\setplotsymbol ({\tiny .})
%
%
%
\plot    2.565 24.871  2.838 24.980
         2.916 25.005
         3.013 25.028
         3.128 25.046
         3.193 25.054
         3.262 25.062
         3.337 25.068
         3.415 25.074
         3.499 25.078
         3.587 25.082
         3.680 25.085
         3.778 25.087
         3.881 25.088
         3.988 25.089
         4.095 25.088
         4.198 25.087
         4.297 25.085
         4.392 25.082
         4.482 25.078
         4.568 25.074
         4.650 25.068
         4.727 25.062
         4.800 25.054
         4.869 25.046
         4.994 25.028
         5.102 25.005
         5.192 24.980
         /
\plot  5.192 24.980  5.520 24.871 /
\linethickness= 0.500pt
\setplotsymbol ({\tiny .})
%
%
%
\plot    2.457 23.669  2.675 23.559
         2.743 23.533
         2.836 23.510
         2.954 23.491
         3.023 23.483
         3.099 23.476
         3.181 23.470
         3.269 23.464
         3.363 23.459
         3.464 23.455
         3.572 23.452
         3.686 23.450
         3.806 23.449
         3.933 23.448
         4.060 23.449
         4.181 23.450
         4.296 23.452
         4.405 23.455
         4.508 23.459
         4.605 23.464
         4.696 23.470
         4.781 23.476
         4.860 23.483
         4.933 23.491
         5.000 23.500
         5.061 23.510
         5.166 23.533
         5.246 23.559
         /
\plot  5.246 23.559  5.520 23.669 /
\linethickness= 0.500pt
\setplotsymbol ({\tiny .})
%
%
%
\plot    1.909 25.307  2.183 25.581
         2.272 25.645
         2.332 25.674
         2.402 25.701
         2.483 25.726
         2.573 25.748
         2.674 25.768
         2.785 25.787
         2.907 25.803
         2.971 25.810
         3.038 25.817
         3.108 25.823
         3.180 25.828
         3.255 25.833
         3.332 25.838
         3.412 25.842
         3.494 25.845
         3.579 25.848
         3.667 25.851
         3.757 25.853
         3.850 25.854
         3.945 25.855
         4.043 25.855
         4.141 25.855
         4.236 25.854
         4.328 25.853
         4.418 25.851
         4.505 25.848
         4.589 25.845
         4.671 25.842
         4.750 25.838
         4.827 25.833
         4.900 25.828
         4.971 25.823
         5.040 25.817
         5.106 25.810
         5.169 25.803
         5.287 25.787
         5.394 25.768
         5.491 25.748
         5.577 25.726
         5.653 25.701
         5.772 25.645
         5.848 25.581
         /
\plot  5.848 25.581  6.068 25.307 /
\linethickness= 0.500pt
\setplotsymbol ({\tiny .})
%
%
%
\plot    1.143 25.199  1.089 25.363
         1.076 25.448
         1.090 25.540
         1.131 25.639
         1.199 25.745
         1.258 25.798
         1.353 25.847
         1.414 25.871
         1.483 25.893
         1.562 25.915
         1.650 25.936
         1.747 25.956
         1.852 25.975
         1.967 25.994
         2.090 26.011
         2.156 26.020
         2.223 26.028
         2.293 26.036
         2.365 26.044
         2.439 26.052
         2.515 26.059
         2.594 26.066
         2.674 26.073
         2.756 26.080
         2.839 26.086
         2.921 26.092
         3.003 26.097
         3.086 26.102
         3.168 26.107
         3.251 26.111
         3.334 26.114
         3.417 26.118
         3.500 26.120
         3.583 26.123
         3.666 26.125
         3.749 26.126
         3.833 26.127
         3.916 26.128
         4.000 26.128
         4.084 26.128
         4.168 26.127
         4.252 26.126
         4.336 26.125
         4.420 26.123
         4.504 26.120
         4.589 26.118
         4.673 26.114
         4.758 26.111
         4.843 26.107
         4.928 26.102
         5.013 26.097
         5.098 26.092
         5.183 26.086
         5.269 26.080
         5.354 26.073
         5.438 26.066
         5.520 26.059
         5.600 26.052
         5.677 26.044
         5.751 26.036
         5.823 26.028
         5.893 26.020
         5.960 26.011
         6.086 25.994
         6.203 25.975
         6.311 25.956
         6.408 25.936
         6.495 25.915
         6.573 25.893
         6.699 25.847
         6.785 25.798
         6.833 25.745
         6.874 25.639
         6.888 25.540
         6.874 25.448
         6.833 25.363
         /
\plot  6.833 25.363  6.723 25.199 /
\linethickness= 0.500pt
\setplotsymbol ({\tiny .})
%
%
%
\plot    0.707 24.651  0.488 24.816
         0.389 24.912
         0.310 25.035
         0.278 25.107
         0.252 25.186
         0.230 25.272
         0.214 25.364
         0.204 25.457
         0.203 25.545
         0.211 25.628
         0.227 25.706
         0.252 25.778
         0.285 25.846
         0.377 25.965
         0.454 26.020
         0.578 26.074
         0.657 26.102
         0.748 26.129
         0.850 26.156
         0.964 26.183
         1.090 26.211
         1.157 26.224
         1.226 26.238
         1.299 26.252
         1.375 26.265
         1.453 26.279
         1.535 26.293
         1.619 26.306
         1.706 26.320
         1.796 26.334
         1.889 26.347
         1.985 26.361
         2.084 26.375
         2.186 26.389
         2.290 26.402
         2.396 26.415
         2.503 26.428
         2.610 26.439
         2.718 26.449
         2.825 26.459
         2.934 26.467
         3.043 26.475
         3.152 26.481
         3.261 26.487
         3.371 26.491
         3.482 26.495
         3.593 26.497
         3.704 26.499
         3.816 26.500
         3.928 26.499
         4.041 26.498
         4.154 26.496
         4.268 26.493
         4.382 26.489
         4.496 26.484
         4.611 26.478
         4.726 26.471
         4.842 26.463
         4.958 26.454
         5.074 26.444
         5.191 26.433
         5.309 26.422
         5.427 26.409
         5.545 26.395
         5.664 26.380
         5.783 26.365
         5.902 26.348
         6.020 26.331
         6.135 26.314
         6.246 26.296
         6.353 26.278
         6.457 26.260
         6.557 26.242
         6.654 26.223
         6.748 26.205
         6.837 26.186
         6.924 26.167
         7.006 26.147
         7.085 26.128
         7.161 26.108
         7.233 26.088
         7.302 26.068
         7.367 26.047
         7.486 26.006
         7.592 25.964
         7.683 25.920
         7.760 25.876
         7.872 25.785
         7.927 25.691
         7.949 25.598
         7.960 25.509
         7.962 25.426
         7.953 25.348
         7.934 25.276
         7.905 25.208
         7.816 25.089
         7.710 24.987
         7.611 24.898
         7.519 24.823
         7.434 24.761
         /
\plot  7.434 24.761  7.271 24.651 /
\linethickness= 0.500pt
\setplotsymbol ({\tiny .})
%
%
%
\plot    0.817 23.887  0.598 23.668
         0.495 23.555
         0.406 23.435
         0.331 23.308
         0.298 23.242
         0.269 23.174
         0.246 23.107
         0.231 23.041
         0.227 22.914
         0.258 22.795
         0.323 22.682
         0.429 22.566
         0.500 22.503
         0.583 22.437
         0.677 22.367
         0.784 22.293
         0.903 22.216
         0.967 22.177
         1.034 22.136
         1.109 22.096
         1.199 22.059
         1.303 22.025
         1.420 21.993
         1.485 21.978
         1.552 21.963
         1.623 21.949
         1.698 21.936
         1.776 21.924
         1.858 21.912
         1.943 21.901
         2.032 21.890
         2.125 21.880
         2.220 21.871
         2.320 21.862
         2.423 21.854
         2.529 21.847
         2.639 21.840
         2.753 21.834
         2.870 21.829
         2.990 21.824
         3.114 21.820
         3.178 21.818
         3.242 21.816
         3.307 21.815
         3.373 21.813
         3.440 21.812
         3.508 21.811
         3.576 21.810
         3.646 21.809
         3.716 21.809
         3.788 21.808
         3.860 21.808
         3.933 21.808
         4.006 21.808
         4.078 21.808
         4.150 21.809
         4.221 21.810
         4.291 21.811
         4.360 21.812
         4.429 21.814
         4.497 21.816
         4.564 21.818
         4.630 21.820
         4.696 21.823
         4.761 21.825
         4.825 21.828
         4.888 21.832
         5.013 21.839
         5.134 21.847
         5.253 21.856
         5.369 21.866
         5.481 21.877
         5.591 21.889
         5.698 21.902
         5.801 21.916
         5.902 21.931
         5.999 21.947
         6.094 21.964
         6.186 21.982
         6.274 22.000
         6.360 22.020
         6.443 22.041
         6.523 22.062
         6.599 22.085
         6.673 22.108
         6.744 22.133
         6.812 22.158
         6.876 22.185
         6.997 22.241
         7.106 22.300
         7.206 22.361
         7.301 22.422
         7.391 22.483
         7.475 22.543
         7.555 22.602
         7.629 22.662
         7.698 22.720
         7.763 22.779
         7.876 22.894
         7.968 23.008
         8.040 23.120
         8.091 23.230
         8.121 23.335
         8.129 23.429
         8.114 23.513
         8.077 23.586
         8.018 23.650
         7.937 23.702
         7.834 23.745
         7.708 23.778
         /
\plot  7.708 23.778  7.161 23.887 /
\linethickness= 0.500pt
\setplotsymbol ({\tiny .})
%
%
%
\plot    1.581 23.338  1.526 23.065
         1.518 22.999
         1.519 22.936
         1.553 22.820
         1.629 22.717
         1.745 22.628
         1.838 22.590
         1.902 22.572
         1.978 22.556
         2.066 22.542
         2.166 22.528
         2.278 22.516
         2.402 22.505
         2.468 22.500
         2.538 22.496
         2.610 22.491
         2.686 22.487
         2.764 22.484
         2.845 22.480
         2.930 22.477
         3.017 22.474
         3.108 22.472
         3.201 22.470
         3.297 22.468
         3.397 22.467
         3.499 22.466
         3.604 22.465
         3.712 22.464
         3.824 22.464
         3.935 22.464
         4.044 22.465
         4.150 22.465
         4.254 22.466
         4.354 22.467
         4.452 22.468
         4.548 22.469
         4.641 22.471
         4.731 22.473
         4.818 22.475
         4.903 22.477
         4.985 22.480
         5.064 22.482
         5.141 22.485
         5.215 22.488
         5.287 22.492
         5.355 22.495
         5.421 22.499
         5.546 22.507
         5.659 22.516
         5.762 22.526
         5.854 22.537
         5.935 22.548
         6.006 22.561
         6.066 22.574
         6.168 22.605
         6.255 22.643
         6.381 22.739
         6.446 22.861
         6.456 22.933
         6.449 23.011
         /
\plot  6.449 23.011  6.394 23.338 /
%
%
\put{\Huge T} [lB] at  1.37 23.997
%
%
\put{\Huge$\mbox{T}^\prime$} [lB] at  6.068 23.997
\linethickness=0pt
\putrectangle corners at  0.169 26.530 and  8.191 21.791
\endpicture}
\end{center}
\caption{A connected sum of two 8-tangles $T$ and $T'$}
\label{fig:question}
\end{figure}

For any positive even integer $p$, choose a fixed set $E$ of $p$
points in $S^2=\partial B$ where $B$ is the unit ball in $\capr^3$.
Define a $p$-tangle as a one-dimensional properly embedded tame
subspace of $B$ whose boundary is $E$.  By connecting the endpoints in
a trivial way, a pair $T$ and $T'$ of $p$-tangles gives rise to a link
$T \# T'$.  For $p=8$ $T \# T'$ is drawn in Figure~\ref{fig:question}.

Is there a finite collection $\ssubp$ of $p$-tangles with the following property?
For any tangle $T$ and integer $d$ such that $d \bigm| |\langle T \# S\rangle |$ for all
$p$-tangles $S \in \ssubp$, we have that $d \bigm| |\langle L\rangle |$ for
any link $L$ which intersects the unit ball $B$ in the one-dimensional
subspace $T$.

If $p=4$ this question is answered by the present
paper in the affirmative: The set $\ssubfour$ can be chosen to consist
of two elements $T_0$ and $T_0^*$ such that $T\#T_0$ is ambient isotopic
to $n(T)$ and $T\#T_0^*$ is ambient isotopic to $d(T)$.  The result
now follows from Theorem~\ref{Thm:9main} version~2.

\begin{figure}
\begin{center}
\leavevmode
\epsfbox{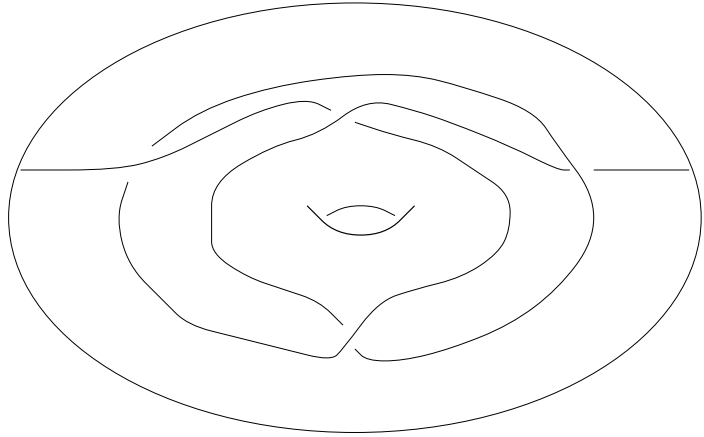}
\end{center} 
\caption{A properly embedded arc in a solid torus}
\label{fig:torus}
\end{figure}

It is natural to ask how a link may intersect subsets of $S^3$ which
are not balls.  For instance, what links can intersect a solid torus
in the arc shown in Figure~\ref{fig:torus}?  A little experimentation
shows that this configuration sits inside a trefoil (determinant $=3$),
a figure-eight knot (determinant $=5$), and a 2,7 torus knot
(determinant $=7$).  Does it sit inside an unknot?

\section{Acknowledgments}

I would like to thank my doctoral thesis advisor Peter B. Shalen for
his patient guidance with the compiling of the thesis which resulted
in this paper.  It was Peter who first recognized the generality of my
methods beyond the case of the link in question being an unknot and
who suggested the ``intersection with cube'' formulation of
embeddability.  In addition I am deeply grateful to Lou Kauffman for
pointing out to me the fraction addition formula of John H. Conway and
for various other suggestions which helped to streamline this paper.
Finally, thanks also to Sandy Rutherford for assistance with \LaTeX.


\begin{thebibliography}{9}

\bibitem{state}

Kauffman, L.:
 State models and the Jones polynomial.
 {\underline{Topology}}, 26(3):395--407, 1987.

\bibitem{jones}

Jones, V.:
 A polynomial invariant for knots via von Neumann algebras.
 {\underline{Bull. Amer. Math. Soc.}}, 12(1):103--112, 1985.

\bibitem{conway}

Conway, J.:
 An enumeration of knots and links and some of their algebraic
  properties.
 In {\underline{Computational Problems in Abstract Algebra}}, ed. J.
  Leech, pages 329--358. New York, Pergammon Press, 1970.

\bibitem{rolfsen}

Rolfsen, D.:
 {\underline {Knots and Links}}.
 Berkeley, CA, Publish or Perish, 1976.

\bibitem{physics}

Kauffman, L.:
 {\underline {Knots and Physics}}.
 Teaneck, NJ, World Scientific, 1991.

\bibitem{rational}

Goldman, J. and Kauffman, L.:
 Rational tangles.
 {\underline{Advances in Applied Mathematics}}, 18:300--332, 1997.

\bibitem{trip}

R.H. Fox:
 A quick trip through knot theory.
 In {\underline{Topology of 3-Manifolds and Related Topics}}, ed. M.K. Fort, pages 120-167. New Jersey, Prentice-Hall, 1961.


\bibitem{silverwms}

Silver, D. and Williams, S.:
 Generalized $n$-colorings of links.
{\underline{Proceedings of 1995 Conference in Knot Theory}},
Banach Center Publications, Vol. 42, 1998.

\bibitem{livingston}

Livingston, C.: {\underline {Knot Theory}}.  Carus Mathematical
 Monographs Vol. 24, Mathematical Association of America, 1993.
 

\end{thebibliography}
\end{document}